\documentclass[journal]{IEEEtran}

\ifCLASSINFOpdf

\else

\fi

\usepackage{mathrsfs}
\usepackage{amsfonts}
\usepackage{amssymb}
\usepackage{amsmath}
\usepackage{cases}
\usepackage{indentfirst}
\usepackage{graphicx}
\usepackage{subfigure}
\usepackage{color}
\usepackage{diagbox}

\newtheorem{thm}{Theorem}
\newtheorem{lem}{Lemma}

\newtheorem{pps}{Proposition}
\newtheorem{dfn}{Definition}
\newtheorem{asm}{Assumption}
\newtheorem{rmk}{Remark}

\newcommand{\tabincell}[2]{\begin{tabular}{@{}#1@{}}#2\end{tabular}}

\begin{document}

\title{Periodic event-triggered networked control systems subject to large transmission delays}

\author{Hao Yu, and Tongwen Chen, {\IEEEmembership{Fellow, IEEE}}
\thanks{This work was supported by NSERC, and an Alberta EDT Major Innovation Fund.}
\thanks{Hao Yu and Tongwen Chen are with the Department of Electrical and Computer Engineering, University of Alberta, Edmonton, AB, T6G 1H9, Canada (e-mail: hy10@ualberta.ca; tchen@ualberta.ca).}}

\markboth{}%
{}

\maketitle

\begin{abstract}
This paper studies periodic event-triggered networked control for nonlinear systems, where the plants and controllers are connected by multiple independent communication channels. Several network-induced imperfections are considered simultaneously, including time-varying inter-sampling times, sensor node scheduling, and especially, large time-varying transmission delays, where the transmitted signal may arrive at the destination node after the next transmission occurs. A new hybrid system approach is provided to model the closed-loop system that contains all communication related behavior. Then, by constructing new storage functions on the system state and updating errors, the relationship between the maximum allowable sampling period (MASP) and maximum allowable delay number in sampling (MADNS) is analyzed, where the latter denotes how many inter-sampling periods can be included in one transmission delay. Moreover, to efficiently reduce unnecessary transmissions, a new dynamic event-triggered control scheme is proposed, where the event-triggering conditions are detected only at aperiodic and asynchronous sampling instants. From emulation-based method, where the controllers are initially designed by ignoring all the network-induced imperfections, sufficient conditions on the dynamic event-triggered control are given to ensure closed-loop input-to-state stability with respect to external disturbances. Moreover, according to different capacities of the communication channels, the corresponding implementation strategies of the designed dynamic event-triggered control schemes are discussed. Finally, two nonlinear examples are simulated to illustrate the feasibility and efficiency of the theoretical results.
\end{abstract}

%
\IEEEpeerreviewmaketitle

\section{INTRODUCTION}
{\IEEEPARstart{I}{n} recent decades, the control and industrial communities have witnessed tremendous development of networked control systems (NCSs), where different components (such as plants and controllers) are physically distributed and connected via (wireless) communication networks~\cite{Nesic2004}. Compared with the traditional dedicated point-to-point and wire-linked control systems, NCSs have many advantages, such as lower costs, reduced weights and volumes, ease of installation and maintenance, and higher reliability, which enable wide applications of NCSs in, e.g., smart grids, wide-area systems, automobiles, and aircrafts \cite{Wangperiodic}.

However, despite the benefits offered by the usage of wireless communication networks, NCSs also suffer some inevitable network-induced imperfections. In this paper, we mainly focus on the following issues: asynchronous and decentralized communication networks, time-varying inter-sampling times and transmission delays, and sensor node scheduling. Due to possibly large scales of NCSs, there may be more than one communication channels, in which, the network and relevant equipment have to operate in an asynchronous and decentralized way. Since networked communication is inherently digital (packet-based), signals cannot be transmitted continuously and instantaneous. Thus, due to the drifting of (independent) local clocks and the varying network conditions, inter-sampling times and transmission delays in all communication channels of NCSs are different and time-varying. Meanwhile, in one communication channel, there could be multiple sensor nodes but only some of them are granted potential access to the network, which results in the design of sensor scheduling protocols. It has been widely recognized that these network-induced imperfections can degrade the control loop performance or even destroy the closed-loop stability. There are several publications on how to understand and compensate the effects of network-induced imperfections, see, e.g., \cite{Heemels2010}, \cite{Dolk Borgers and Heemels}, and some overview papers \cite{Geoverview}, \cite{Liuoverview}.

A fundamental and important aspect in the study of NCSs is on the transmission delays, because, different from the other network-induced issues, transmission delays would influence the real-time capability of system operations. Based on the relationship with inter-sampling periods, the transmission delays are often classified into two cases: the small-delay case and large-delay case \cite{Liuoverview}. In the former, the delay of one transmission has to be no larger than the corresponding inter-sampling period; otherwise, it is in the large-delay case. In linear NCSs, there are several well-developed technical tools to deal with both cases, such as the discrete-time modeling approach \cite{Zhang2001}, time-delay approach \cite{Fridman2004}, and mixed interval-integral approach \cite{Xiao2020}, which directly depended on dynamics and solution structures of linear systems. As a result, it is hard to extend these methods to nonlinear systems that have non-globally Lipschitz dynamics, except for some applications in parabolic partial differential equation systems, see, e.g., \cite{FridmanPDE}. Especially, when considering time-varying transmission delays and sensor node scheduling simultaneously, till now only the time-delay approach is applicable for linear NCSs in the large-delay case \cite{Liuoverview}. For general nonlinear NCSs, \cite{Naghshtabrizi} developed an hybrid/impulsive system approach based on the emulation-based method \cite{Walsh2001}, where the controllers are initially designed by assuming perfect point-to-point links; then, a resultant discontinuous Lyapunov function (functional) \cite{Naghshtabrizi2010} was provided to characterize the effects of several network-induced imperfections on system stability and performance \cite{Heemels2010}, \cite{DolkAuto2017}, \cite{Abdelrahim2019}. In all the above studies, only small transmission delays were considered for general nonlinear NCSs.  To the best of the authors' knowledge, there is no systematical modelling framework for general nonlinear NCSs in the large-delay case, which is one of the motivation of this study.

Besides the network-induced issues, another potential challenge in NCSs is resource constraints, such as, limited communication bandwidth and restricted energy of onboard batteries. To avoid overusing the network, an attractive solution, event-triggered control (ETC), has been proposed in the last two decades \cite{Tabuada}.  In ETC, the executions of transmission tasks are decided by some online events, which is evaluated by the so-called event-triggering conditions, rather than the elapse of some offline designed periods. Thus, by constructing a closed loop from real-time system behavior to transmission decisions, ETC can strike a more desirable balance between the system performance and resource consumption.

In early studies of ETC, the events need to be evaluated continuously, which somewhat contradicts the digital nature of NCSs. Thus, in \cite{Heemelsperiodic}, a new scheme, periodic event-triggered control (PETC), was designed for linear NCSs, where the events were only detected at some discrete time instants. In the large-delay case, \cite{Xiao2020} studied PETC for linear NCSs with multiple asynchronous transmission channels. For nonlinear NCSs, an over-approximation technique was proposed in \cite{Borgersnonlinear} for converting continuous event-triggered controllers to periodic ones. In \cite{Wangperiodic}, the hybrid system approach was applied to the design of PETC with considerations of time-varying inter-sampling times and sensor node scheduling. Note that all the schemes in the aforementioned publications on PETC belong to the type of static ETC, where only the current sampled signals are involved in the event-triggering conditions. It has been illustrated that the transmission performance in ETC can be improved by designing an auxiliary dynamical system to record the historical online information, which results in the introduction of dynamic ETC \cite{Girard}. Although more and more studies focus on dynamic ETC in recent years \cite{Dolk Borgers and Heemels}, \cite{DolkAuto2017}, \cite{YuHaoandChen}, almost all of them required continuous detection of events \cite{Ge2021}. In \cite{BorgersRiccati} and \cite{Fu2020}, a Riccati-based design method of dynamic PETC was proposed for linear NCSs subject to small (constant) transmission delays. However, due to the dependence of the considered Riccati equation on linear dynamics, their results are hardly possible to be applied to nonlinear NCSs, which gives another motivation of this paper.

Based on the observations above, this paper studies (dynamic) PETC for nonlinear NCSs subject to several network-induced imperfections, especially including the large time-varying transmission delays. The main contributions of this paper are summarized as follows.

First, a new hybrid system approach is provided to model nonlinear NCSs with multiple independent communication channels, which suffer simultaneously time-varying inter-sampling times, sensor node scheduling, and large transmission delays. Then, by constructing new storage functions on the system state and updating errors (differences between the current and most recently updated signals), the relationship between the maximum allowable sampling period (MASP) \cite{Nesic Teel and Carnevale} and maximum allowable delay number in sampling (MADNS) is analyzed, where the latter denotes how many inter-sampling periods can be included in one transmission delay. The proposed modelling and analysis approach includes the one in \cite{Dolk Borgers and Heemels} for small transmission delays as a special case.

Second, to efficiently reduce unnecessary transmissions, a new dynamic PETC scheme is proposed, where the evaluation of event-triggering conditions involves an auxiliary discontinuous variable that records the effects of historical online information by its flow (described by differential equations) and jump (described by difference equations) behavior. From some assumptions provided by the emulation-based method, sufficient conditions on the design of dynamic PETC are given to ensure input-to-state stability of closed-loop systems with respect to external disturbances. Furthermore, according to different capacities of communication channels, corresponding implementation strategies are discussed following the principle that event triggers (ETs), the hardware to realize the detection of events, can utilize only the sampled local information in a way decided by the assumed capacities of communication channels. Furthermore, it is showed that the designed dynamic PETC can lead to better transmission performance than the static one in \cite{Wangperiodic}.

\indent The rest of the paper is organized as follows. After reviewing the basic definitions and notations in Section II, the problem of PETC for nonlinear NCSs subject to large transmission delays is formulated in Section III. The main results of this paper are given in Section IV. Nonlinear examples are given to illustrate the feasibility of the theoretical results in Section V.  The conclusions of this paper are drawn in Section VI.

\section{Notations}
Let $\mathbb{R}$ ($\mathbb{Z}$) be the set of real numbers (integers). $\mathbb{R}_{(\cdot)}$ ($\mathbb{Z}_{(\cdot)}$) denotes all the reals (integers) that satisfy the relationship in $(\cdot)$. The absolute value of a scalar $s\in\mathbb{R}$ is denoted by $|s|$, and Euclidean norm of a vector $x \in \mathbb{R}^n$ is denoted by $\left\| {x} \right\|$. The Euclidean induced matrix norm of $A\in \mathbb{R}^{n \times m}$ is denoted by $\left\| {A} \right\|$. The transpose of a matrix $A \in \mathbb{R}^{n \times m}$ is denoted by $A^{\rm{T}}$. Define $\bar{N}:=\{1,\dots,N\}$ with some $N\in\mathbb{Z}_{\ge1}$. Denote a $N$-dimension column vector with all entries equal to one by $\mathbf{1}_N$. Let $I_n$ ($0_{n\times m}$) be the identity (zero) matrix with dimension $n$ ($n\times$ m), and sometimes the subscript will be omitted if there is no confusion. ${\rm{diag}}(\cdots)$ denotes a diagonal matrix, with diagonal entries listed. Let $|\Omega|$ be the cardinality of the set $\Omega$. For two sets $\Omega_1,\Omega_2\in\mathbb{R}^n$, define $\Omega_1\backslash\Omega_2:=\{x\in\mathbb{R}^n|x\in\Omega_1,x\notin\Omega_2\}$. Given a set $\Omega\subset\mathbb{R}^n$ and a vector $x\in\mathbb{R}^n$, the distance of $x$ to $\Omega$ is defined as $\left\|x\right\|_{\Omega}:=\inf_{y\in\Omega}\left\|x-y\right\|$. The Clarke derivative \cite{Clarke} is defined as follows: for a locally Lipschitz function $U : \mathbb{R}^n \to \mathbb{R}$ and a vector $v\in \mathbb{R}^n$,
\[U^{\circ}(x,v) := \limsup _{h \to 0^+,y \to x}\frac{{U(y + hv) - U(y)}}{h}.\]
If the function $U(\cdot)$ is differential at the point $x$, the Clark derivative $U^{\circ}(x,v)$ reduces to the standard directional derivative $\left\langle\nabla U(x),v\right\rangle$, where $\nabla U(\cdot)$ is the classical gradient.

A continuous function $\gamma: \mathbb{R}_{\ge0} \to \mathbb{R}_{\ge0}$ is called a $\mathcal{K}$--function (denoted by $\gamma\in\mathcal{K}$) if it is continuous, strictly increasing and satisfy $\gamma(0)=0$; $\gamma$ is a $\mathcal{K}_\infty$--function (denoted by $\gamma\in\mathcal{K}_{\infty}$), if $\gamma\in\mathcal{K}$ and $\gamma(s)\to \infty$ as $s\to \infty$. A continuous function $\beta: \mathbb{R}_{\ge 0} \times \mathbb{R}_{\ge 0} \to \mathbb{R}_{\ge 0}$ is a $\mathcal{KL}$--function (denoted by $\beta\in\mathcal{KL}$), if it satisfies: (i) for each $t\ge0$, $\beta(\cdot,t)\in\mathcal{K}$, and (ii) for each $s \ge 0$, $\beta(s,\cdot)$ is nonincreasing and $\lim_{t\to\infty}\beta(s,t)=0$. A function $f:\mathbb{R}^n\to\mathbb{R}$ is locally Lipschitz if for any compact set $\mathcal{S}\subset \mathbb{R}^n$, there exists a constant $L(\mathcal{S})$ such that $\left|f(x)-f(y)\right|\le L\left\|x-y\right\|$ for all $x,y\in\mathcal{S}$; in addition, if $\mathcal{S}=\mathbb{R}^n$, $f$ is globally Lipschitz.  For a function $\omega:\mathbb{R}_{\ge0}\to\mathbb{R}^n$, $\omega(r^{+})$ denotes the limit from above at time $r\in\mathbb{R}_{\ge0}$, i.e., $\omega(r^{+})=\lim_{t\searrow r} \omega(t)$.

Closed-loop systems in this paper will be modelled as a nonlinear hybrid system of the following form:
\begin{equation}
  \mathcal{H}:\begin{cases}
    \dot{q}=F(q,w), &q\in\mathcal{C},\\
    q^+\in G(q), &q\in\mathcal{D},
  \end{cases}
\end{equation}
where $q\in\mathcal{X}\subset\mathbb{R}^{n_q}$ is a state vector and $w\in\mathbb{R}^{n_w}$ is an external input. $F$ is a continuous flow map and $G$ is a outer semi-continuous\footnote{See \cite{Goebelbook} for the definition of outer semi-continuity of set-valued maps.} set-valued jump map. The flow set $\mathcal{C}$ and jump set $\mathcal{D}$ are closed, and the set of initial states, $\mathcal{X}_0$, satisfies $\mathcal{X}_0\subset \mathcal{C}\cup \mathcal{D}\subset\mathcal{X}$. A solution to the hybrid system $\mathcal{H}$ is defined on the hybrid time domain, denoted by ${\text{dom }} q$, where the element is expressed as $(t,\bar{j})$ with $t\in\mathbb{R}_{\ge0}$ and $\bar{j}\in\mathbb{Z}_{\ge0}$ recording, respectively, the elapse of time and number of jumps. The hybrid time domain, ${\text{dom }} q$, will not be often mentioned explicitly if there is no confusion from the context. To save space, we omit the mathematical definitions on hybrid time domain and other notations in hybrid systems, and refer the readers to \cite{Goebelbook}.

\section{Problem formulation}
In this section, we first introduce the configuration of nonlinear NCSs with multiple independent communication channels. Then, in each communication channel, the considered  network-induced issues (time-varying inter-sampling times, sensor node scheduling and large transmission delays) as well as dynamic PETC are introduced. Furthermore, by providing some properties on the evolution of updating errors in the large-delay case, a hybrid system model is provided for closed-loop systems. Finally, the objective of this paper is elaborated.
\subsection{Networked Control Systems}
\begin{figure}[!hbtp]
\centering
\includegraphics[width=8.5cm]{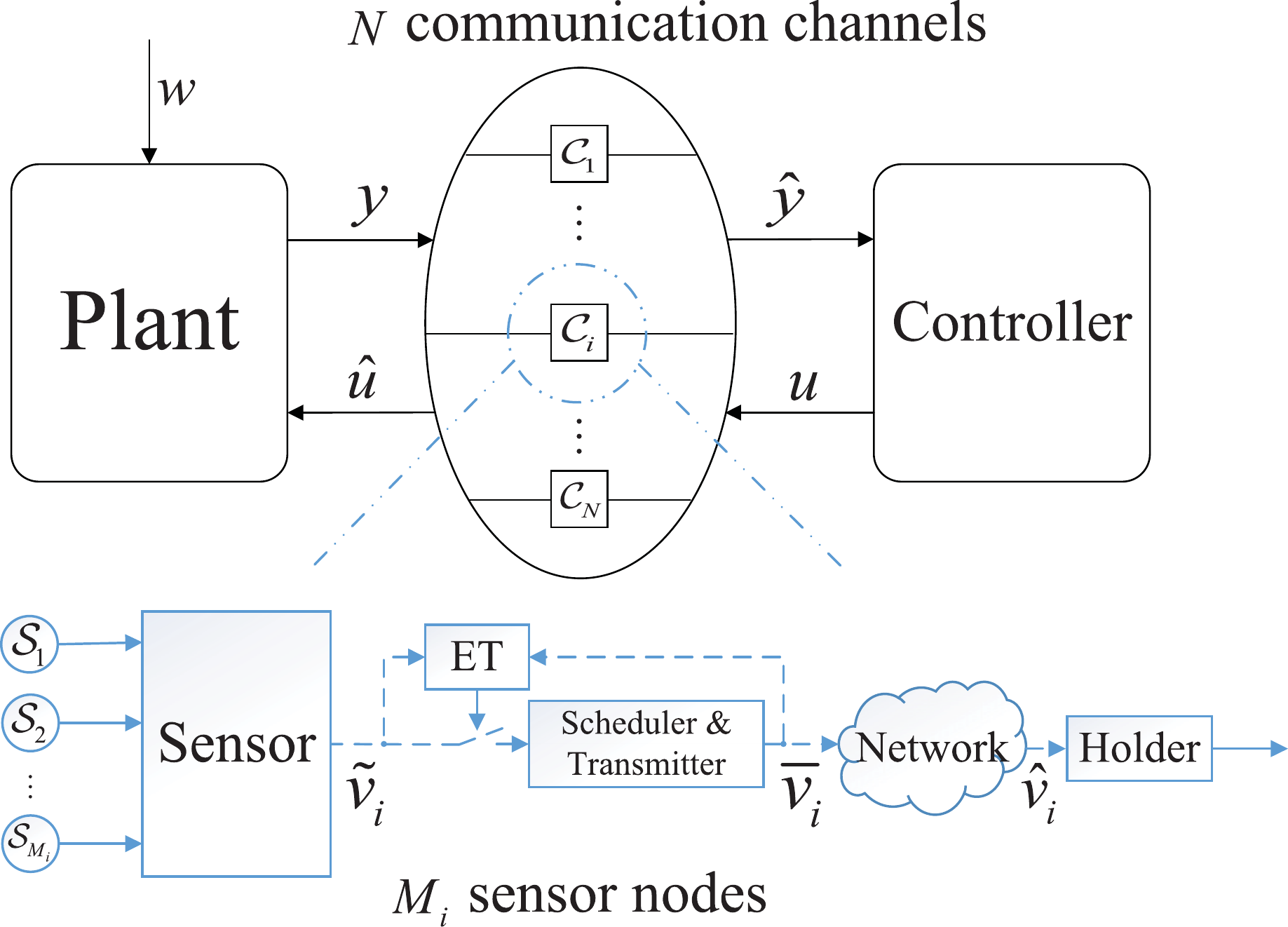}
\caption{Configuration of event-triggered networked control systems, where the signals with solid lines are transmitted continuously and the ones with dashed lines are accessible only at discrete instants.}
\label{configuration}
\end{figure}

Referring to \cite{Dolk Borgers and Heemels}, we consider an event-triggered NCS in Fig.~\ref{configuration}, where the nonlinear time-invariant plant is described by
\begin{equation}\label{plant}
   \dot{x}_p=f_p(x_p,\hat{u},w), {\text{\quad}} y=g_p(x_p),
\end{equation}
where $x_p\in\mathbb{R}^{n_p}, y_p\in\mathbb{R}^{n_y}$ and $w\in\mathbb{R}^{n_w}$ stand for, respectively, the state, measurable output, and bounded external disturbance. The left-continuous and piece-wise constant vector $\hat{u}\in\mathbb{R}^{n_u}$ is the received control signal (in a zero-order-hold manner), with the original version $u\in\mathbb{R}^{n_u}$ generated by the following output-feedback controller:
\begin{equation}\label{controller}
   \dot{x}_c=f_c(x_c,\hat{y}), {\text{\quad}} u=g_c(x_c),
\end{equation}
with the controller state $x_c\in\mathbb{R}^{n_c}$, controller output (original input) $u$, and the received output $\hat{y}$ that is left-continuous and piece-wise constant. The continuous functions $f_p$ and $f_c$ and the continuously differentiable functions $g_p$ and $g_c$ are zero at zero. The controller in (\ref{controller}) has been designed in an emulation-based manner, rending some desirable properties (which will be specified later) for the system without network-induced constraints (i.e., $\hat{y}(t)=y(t),\hat{u}(t)=u(t)$ for all $t\in\mathbb{R}_{\ge0}$).

From Fig. \ref{configuration}, the transmitted signals $v:=[y^{\rm{T}},u^{\rm{T}}]^{\rm{T}}\in\mathbb{R}^{n_v}$, with $n_v:=n_y+n_u$, are communicated through $N$ independent channels $\mathcal{C}_i, i\in\bar{N}$, which leads to the received version $\hat{v}:=[\hat{y}^{\rm{T}},\hat{u}^{\rm{T}}]^{\rm{T}}\in\mathbb{R}^{n_v}$. To facilitate the description in the rest, we partition $v$ ($\hat{v}$) into $v=[v_1^{\rm{T}},\dots,v_N^{\rm{T}}]^{\rm{T}}$ ($\hat{v}=[\hat{v}_1^{\rm{T}},\dots,\hat{v}_N^{\rm{T}}]^{\rm{T}}$) where we assume that $v_i\in\mathbb{R}^{n_{v,i}},i\in\bar{N}$, is transmitted over the communication channel $\mathcal{C}_i$ (possibly after reordering the channel indices). Note that $n_{v,i}$ is not necessarily equal to $1$. In addition, one has $v_i=g_{v_i}(x)$ for some continuously differentiable function $g_{v_i}$ with $g_{v_i}(0)=0$ and $i\in\bar{N}$, where $x:=[x_p^{\rm{T}},x_c^{\rm{T}}]^{\rm{T}}\in\mathbb{R}^{n_x}$ and $n_x=n_p+n_c$.

\begin{rmk}
  \label{rmk of zero controller}
  {Similar to \cite{Dolk Borgers and Heemels}, the considered controller in (\ref{controller}) can include the static case, that is, $n_c=0$ and $u=g_c(\hat{y})$. Since $u$ has been in a sampled-data form and only static calculation on $\hat{y}$ is required, one can suppose that there is no extra digital channels from the controller to the plant, namely, $v=y$ and $u=\hat{u}$.}
\end{rmk}

\subsection{Large Network-induced Delays and Scheduling Protocols}
For the communication channel $\mathcal{C}_i, i\in\bar{N},$ in Fig. \ref{configuration}, the signal $v_i$ is sampled at the sampling instant $s_j^i$, with $s_0^i=0$ for all $i\in\bar{N}$. Assume the following constraints of sampling instants:
\begin{equation}
  \nonumber
  T_m^i\le \tau_j^i \le T_M^i,
\end{equation}
with the inter-sampling time $\tau_j^i:=s_{j+1}^i-s_j^i$, the upper bound $T_M^i>0$ (also known as MASP), and lower bound $T_m^i\in(0,T_M^i]$ of inter-sampling times. Due to the independence of communication channels and drifting of local clocks, the sampling time sequence $\{s_j^i\}_{j\in\mathbb{Z}_{\ge0}}$ is not necessarily periodic or synchronized with others, although we still call the to-be-studied ETC scheme as PETC. For avoiding Zeno behavior \cite{Borgers}, $T_m^i\in(0,T_M^i],i\in\bar{N},$ can be selected arbitrarily small in theory and decided by hardware constraints in reality. The upper bounds $T_M^i,i\in\bar{N},$ essentially affect closed-loop stability and will be designed later.

At each sampling instant, the ET in the communication channel $\mathcal{C}_i,i\in\bar{N}$, in Fig. \ref{configuration} will check the event-triggering condition to decide whether to transmit the current signal $v_i$ through the network. Thus, denote by $\{t_k^i\}_{k\in\mathbb{Z}_{\ge0}}\subset\{s_j^i\}_{j\in\mathbb{Z}_{\ge0}}$ the transmission time sequence of $\mathcal{C}_i,i\in\bar{N}$ with $t_0^i=s_0^i$. A general form of event-triggering condition is described as
\begin{equation}
  \label{ET general}
  t_{k+1}^i=\min\{t>t_k^i|t\in\{s_j^i\}_{j\in\mathbb{Z}_{\ge0}}, g_s^i(o_i(t),\eta_i(t))<0\},
\end{equation}
where the local information vector $o_i\in\mathbb{R}^{n_{o,i}}$ will be specified later and $g_s:\mathbb{R}^{n_{o,i}}\times \mathbb{R}_{\ge0}\to\mathbb{R}$ is the to-be-designed triggering function. The left-continuous auxiliary variable $\eta_i\in\mathbb{R}_{\ge0}$ follows the following dynamics
\begin{equation}
  \label{dynamics of etai}
  \begin{cases}
    \dot{\eta}_i=f_{\eta}^i(o_i(t),\eta_i(t)),  & t\in(s_j^i,s_{j+1}^i],\\
    \eta_i(t^+)=g_t^i(o_i(t),\eta_i(t)), & t\in\{t_k^i\}_{k\in\mathbb{Z}_{\ge0}},\\
    \eta_i(t^+)=g_s^i(o_i(t),\eta_i(t)), & t\in\{s_j^i\}_{j\in\mathbb{Z}_{\ge0}}\backslash \{t_k^i\}_{k\in\mathbb{Z}_{\ge0}},\\
  \end{cases}
\end{equation}
with initial state $\eta_i(0)=0$, where the functions $f_{\eta}^i:\mathbb{R}^{n_{o,i}}\times \mathbb{R}_{\ge0}\to\mathbb{R}$ and $g_t^i:\mathbb{R}^{n_{o,i}}\times \mathbb{R}_{\ge0}\to\mathbb{R}$ satisfy that $f_{\eta}^i(\cdot, 0)$ and $g_t^i(\cdot,\eta_i)$ for all $\eta_i\ge0$ are non-negative, which combined with (\ref{ET general}) ensures the non-negativeness of $\eta_i,i\in\bar{N}$.

Due to the communication over networks, the update of $\hat{v}_i$ corresponding to the transmission time $t_k^i$ for any $k\in\mathbb{Z}_{\ge0}$ and $i\in\bar{N}$ suffers a network-induced delay $d_k^i\in\mathbb{R}_{\ge0}$, which results in the arrival time $f_k^i=t_k^i+d_k^i$. Different from \cite{Heemels2010}, in this paper, we consider the following large-delay case, where the arrival time $f_k^i$ can be larger than the next sampling instant $t_k^i+\tau_j^i$, or even lager than the next transmission instant $t_{k+1}^i$ with some $j\in\mathbb{Z}_{\ge0}$ satisfying $t_k^i=s_{j}^i$. Specifically, we introduce the following assumption.

\begin{asm}
  \label{asm of large delay}
  For the communication channel $\mathcal{C}_i,i\in\bar{N}$, there exists an integer $D_i\in\mathbb{Z}_{\ge0}$, called maximum allowable delay number in sampling (MADNS), such that $0\le d_k^i \le \sum_{n=0}^{D_i}\tau_{j+n}^i$ for all $k,j\in\mathbb{Z}_{\ge0}$ satisfying $t_k^i=s_{j}^i$.
\end{asm}

Assumption \ref{asm of large delay} includes the small-delay case as a special scenario by choosing $D_i=0$. Due to the independency of communication channels, it is allowable to assume $D_i\neq D_j$ for $i,j\in\bar{N}$ and $i\neq j$. Moreover, {similar to Assumption 1 in \cite{Xiao2020}}, it is suppose that there is no disorder in updates, that is, in each $\mathcal{C}_i,i\in\bar{N}$, earlier transmitted signal $v_i$ arrives at the destination node (i.e., the holder at the other side of the network in Fig. \ref{configuration}) also earlier.

\begin{rmk}\label{rmk of delay assumption}
  In Assumption \ref{asm of large delay}, the large delay is counted by the numbers in sampling instead of the absolute values of $\max\{f_{k+1}-f_k\}_{k\in\mathbb{Z}_{\ge0}}$. Thus, a system in the small-delay case with a large inter-sampling time could yield a large value of delay than that in the large-delay case but with a smaller inter-sampling time. This is reasonable because the maximum allowable time of transmission delays is mainly decided by the system dynamics rather than the sampling frequency. Consequently, the value of MADNS should decrease as the sensor samples slower.
\end{rmk}

Moreover, in Fig. \ref{configuration}, we consider $M_i\in\{1,\dots, n_{v,i}\}$ sensor nodes in $\mathcal{C}_i,i\in\bar{N}$. Due to the limited capacity of networks, at each possible transmission instant, the transmitter can only send the information collected by parts of the nodes. Thus, the transmission of $v_i$ involves the so-called scheduling protocols, as introduced in \cite{Nesic2004}, to determine which nodes are granted potentially access to the network. In detail, we have that, for each communication channel $i\in\bar{N}$,
\begin{equation}
  \label{scheduling protocol}
  \begin{split}
    \bar{v}_i\left(t_k^{i+}\right)=v_i(t_k^i)+h_{v_i}(k,\bar{e}_i(t_k^{i})),
  \end{split}
\end{equation}
where the piece-wise constant and left-continuous vector $\bar{v}_i(t)=\bar{v}_i(t_k^i), t\in(t_k^i,t_{k+1}^i],$ records the scheduled signal, and the transmission error is defined as $\bar{e}_i(t):=\bar{v}_i(t)-v_i(t)$. The function $h_{v_i}:\mathbb{Z}_{\ge0} \times \mathbb{R}^{n_{v,i}}\to\mathbb{R}^{n_{v,i}}$ is an update function decided by scheduling protocols, such as the popular sampled-data (SD), round-robin, and try-once-discard protocols \cite{Nesic2004}. Consequently, the update of $\hat{v}_i$ at the destination node of $\mathcal{C}_i$, satisfies
\begin{equation}
  \label{update of hat vi}
  \hat{v}_i\left(f_k^{i+}\right)=\bar{v}_i\left(t_k^{i+}\right)=v_i(t_k^{i})+h_{v_i}(k,\bar{e}(t_k^{i})).
\end{equation}
Since $f_0^i$ might be larger than $0$, we assume that the initial value of $\hat{v}_i(0)$ is known both for the ET and holder. This can be trivially achieved by selecting zero initial values for each communication channel $\mathcal{C}_i,i\in\bar{N}$.

Moreover, define the left-continuous vector $\tilde{v}_i(t),i\in\bar{N}$, to record the latest sampled signal:
\begin{equation}
  \nonumber
  \tilde{v}_i(t)=v_i(s_j^i), t\in(s_j^i,s_{j+1}^i].
\end{equation}
Note that $\tilde{v}_i$ does not involve $h_{v_i}$ since it will be used by the ET instead of the transmitter. An illustration of the time sequences $\{s_j^i\}_{j\in\mathbb{Z}_{\ge0}}$, $\{t_k^i\}_{k\in\mathbb{Z}_{\ge0}}$, $\{f_k^i\}_{k\in\mathbb{Z}_{\ge0}}$ and different versions of $v_i$ is showed in Fig. \ref{TS}.

{Therefore, the interest of this paper is to design the decentralized PETC in (\ref{ET general}), namely, constant $T_M^i$ and functions $f_\eta^i,g_t^i,g_s^i$, for all communication channels $i\in\bar{N}$ to ensure the input-to-state stability, whose formal and exact definition will be given in Section~III--D.}

\begin{figure*}[!hbtp]
\centering
\includegraphics[width=13cm]{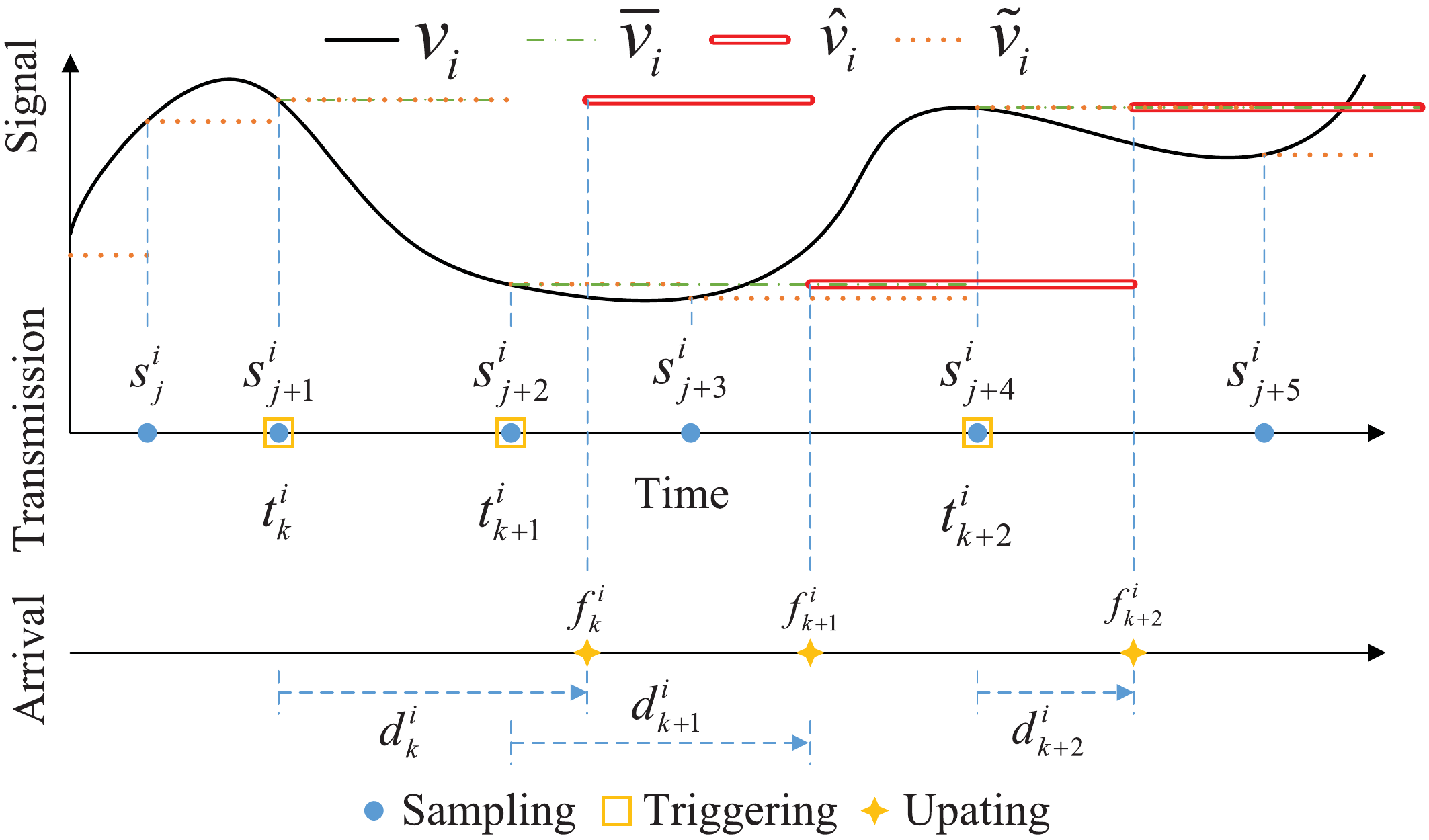}
\caption{Illustration for different time sequences and signals in $\mathcal{C}_i,i\in\bar{N}$ with $D_i=1$ and $M_i=1$ (resulting in $h_{v_i}=0$).}
\label{TS}
\end{figure*}

\subsection{Evolution of Updating Errors}
Due to the large-delay assumption and scheduling protocols, the relationship between the transmission error $\bar{e}_i$ and the updating error $e_i:=\hat{v}_i-v_i$ becomes more complex than that in the small-delay case. To study the evolution of $e_i,i\in\bar{N}$, we introduce the following left-continuous variables in Table~\ref{table1} for each communication channel $\mathcal{C}_i$, which leads to $\bar{k}_i(t_k^i)=k$ and $\tilde{k}_i(f_k^i)=k$ for all $i\in\bar{N}$. Moreover, it is worth noting that $l_i(s_j^i)\le D_i$, since, from Assumption \ref{asm of large delay}, the transmission before $s_{j-D_i}^i$ must have been updated and $l_i(t)$ is left-continuous.
\renewcommand\arraystretch{1.5}
\begin{table}[!hbtp]
\caption{{Definitions of  variables}.}
\begin{center}
\begin{tabular}{c|l}
\hline
 Variables & Description  \\
\hline
 $\bar{k}_i(t)\in\mathbb{Z}_{\ge 0}$& total transmission number in $\mathcal{C}_i$ up to $t$ \\
 $\tilde{k}_i(t)\in\mathbb{Z}_{\ge 0}$& total updating number in $\mathcal{C}_i$ up to $t$ \\
 $l_i(t)\in\{0,\dots,D_i+1\}$  & \tabincell{l}{difference between transmission and updating\\ numbers: $\bar{k}_i(t)-\tilde{k}_i(t)$}\\
 $\hat{m}_i(t)\in\{-1,1\}$   & \tabincell{l}{the next action is transmission ($\hat{m}=1$) or\\ updating ($\hat{m}=-1$)}\\
 \hline
\end{tabular}
\end{center}
\label{table1}
\end{table}

From Table \ref{table1}, the most recently received signal $\hat{v}_i,i\in\bar{N}$, satisfies that for $t=f_k^i$, one has $\hat{m}_i(t)=-1$, $\tilde{k}_i(t^+)=\tilde{k}_i(t)+1$ and
\begin{equation}
  \label{update of hat v}
  \begin{split}
    \hat{v}_i\left(t^+\right)=& \bar{v}_i\left(t_{\tilde{k}_i(t)}^{i+}\right)=\bar{v}_i\left(t_{\bar{k}_i(t)-l_i(t)}^{i+}\right).
  \end{split}
\end{equation}
Consequently, applying (\ref{update of hat v}) to the transmission error $e_i$ yields that for $t=f_k^i$,
\begin{equation}
  \label{evolution of ei}
  \begin{split}
      e_i(t^+)=&\hat{v}_i(t^+)-v_i(t)+\hat{v}_i(t)-\hat{v}_i(t)\\
      =&\hat{v}_i(f_{k-1}^{i+})-v_i(t)+\hat{v}_i(f_k^{i+})-\hat{v}_i(f_{k-1}^{i+})\\
      =&e_i(t)+\bar{v}_i\left(t_{\tilde{k}_i(t)}^{i+}\right)-\bar{v}_i\left(t_{\tilde{k}_i(t)-1}^{i+}\right)\\
      =&e_i(t)+\bar{v}_i\left(t_{\bar{k}_i(t)-l_i(t)}^{i+}\right)-\bar{v}_i\left(t_{\bar{k}_i(t)-l_i(t)-1}^{i+}\right).
  \end{split}
\end{equation}
Based on (\ref{evolution of ei}), for all $i\in\bar{N}$, we introduce the following memory vector $\theta_i(t)=\left[\theta_{i,1}^{\rm{T}}(t),\dots,\theta_{i,D_i+1}^{\rm{T}}(t) \right]^{\rm{T}}\in\mathbb{R}^{n_{\theta,i}}$ with $n_{\theta,i}:=n_{v,i}(D_i+1)$, satisfying
\begin{equation}
\label{def of theta}
  \theta_{i,j}(t)=\bar{v}_i\left(t_{\bar{k}_i(t)-l_i(t)+j-1}^{i+}\right)-\bar{v}_i\left(t_{\bar{k}_i(t)-l_i(t)+j-2}^{i+}\right),
\end{equation}
for $j=1,\dots, l_i(t)$; otherwise, $\theta_{i,j}(t)=0$, which leads to
\begin{equation}
  \label{update of ei}
  e_i(f_k^{i+})=e_i(f_k^i)+\theta_{i,1}(f_k^i).
\end{equation}
Moreover, for $\theta_i(t),i\in\bar{N}$, we propose the following property.
\begin{pps}\label{pps of thetai}
  The following three statements hold for $\theta_i(t),i\in\bar{N}$:
  \begin{enumerate}
    \item if $t\in\{s_j^i\}_{j\in\mathbb{Z}_{\ge0}}\backslash\{t_k^i\}_{k\in\mathbb{Z}_{\ge0}}$, one has $\theta_{i}(t^+)=\theta_{i}(t)$;
    \item if $t=f_k^i$, one has $\theta_{i,j}(t^+)=\theta_{i,j+1}(t)$ for $j=1,\dots, D_i$ and $\theta_{i,D_i+1}(t)=0$;
    \item if $t=t_k^i$, one has $\theta_{i,j}(t^+)=\theta_{i,j}(t)$ for $j\in\{1,\dots,D_i+1\}\backslash\{l_i(t)+1\}$ and
    \[\begin{split}
      \theta_{i,l_i(t)+1}(t^+)=&h_{v_i}(\bar{k}(t),\bar{e}_i(t_{\bar{k}(t)}^{i}))-e_i(t)-\sum_{j=1}^{l_i(t)}\theta_{i,j}(t),\\
      =&h_{v_i}(\bar{k}(t),\bar{e}_i(t_{\bar{k}(t)}^{i}))-\bar{e}_i(t_{\bar{k}(t)}^{i}).
    \end{split}\]
  \end{enumerate}
\end{pps}
\begin{IEEEproof}
  See Appendix B.
\end{IEEEproof}

\begin{rmk}\label{rmk of pps of thetai}
  Since there is no disorder in updates, the first block in $\theta_i$ always stores the information that will be invoked in the next updates, namely (\ref{update of ei}). Thus, after an updating instant, the used information $\theta_{i,1}$ will be discarded (Item 2)). At a sampling instant that is a transmission instant as well, a new vector $\theta_{i,l+1}$ should be stored following the equality in Item 3) for future use; otherwise, no information is memorized (Item 1)).
\end{rmk}

Consequently, for each communication channel $\mathcal{C}_i,i\in\bar{N}$, the local information of its ET is
\[o_i:=\left[\tilde{v}_i,\bar{k}_i,\bar{v}_i,\theta_i^{\rm{T}},\eta_i,l_i,e_i\right]^{\rm{T}}.\]
Note that the implementation of some information in $o_i$ depends on the capacity of equipment. For example, $\tilde{k}_i$ and $l_i$ are unknown by the ET without delay-free acknowledgment mechanism \cite{DolkAuto2017}.

\subsection{Closed-loop System Models}
Define the following augmented states:
\[\begin{split}
  e:=&[e_1^{\rm{T}},\dots,e_N^{\rm{T}}]^{\rm{T}},{\text{ }}\theta:=[\theta_1^{\rm{T}},\dots,\theta_N^{\rm{T}}]^{\rm{T}}, {\text{ }}\hat{\tau}:=[\hat{\tau}_1,\dots,\hat{\tau}_N]^{\rm{T}},\\
  {\text{ }}\bar{k}:=&[\bar{k}_1,\dots,{\text{ }}\bar{k}_N]^{\rm{T}},\tilde{k}:=[\tilde{k}_1,\dots,\tilde{k}_N]^{\rm{T}},\eta:=[\eta_1,\dots,\eta_N]^{\rm{T}},\\
   l:=&[l_1,\dots,l_N]^{\rm{T}},{\text{ }}\hat{m}:=[\hat{m}_1,\dots,\hat{m}_N]^{\rm{T}},{\text{ }} \tilde{v}=[\tilde{v}_1^{\rm{T}},\dots,\tilde{v}_N^{\rm{T}}]^{\rm{T}},
\end{split}\]
where $\hat{\tau}_i(t):=t-s_j^i$ for $t\in(s_j^i,s_{j+1}^i]$ records the time elapsed since the last sampling instants in the communication channel $\mathcal{C}_i,i\in\bar{N}$.  Then, the closed-loop system, with the state $q:=\left[x^{\rm{T}},e^{\rm{T}},\theta^{\rm{T}},\tilde{v}^{\rm{T}},\hat{\tau}^{\rm{T}},\bar{k}^{\rm{T}},\tilde{k}^{\rm{T}},l^{\rm{T}},\hat{m}^{\rm{T}},\eta^{\rm{T}}\right]^{\rm{T}}\in\mathcal{X}$ and
\[\begin{split}
 \mathcal{X}:=&\{q\in\mathcal{\chi}_{\rm{all}}|l_i\in\{1,\dots,D_i\},\theta_{i,j}=0,\\
 &{\text{\quad}}{\text{\quad}}{\text{\quad}}{\text{\quad}}{\text{for }} j=l_i+1,\dots,D_i+1,i\in\bar{N}\},\\
 \mathcal{X}_{\rm{all}}:=&\mathbb{R}^{n_x}\times\mathbb{R}^{n_v}\times \mathbb{R}^{n_{\theta}} \times \mathbb{R}^{n_v}\times \mathbb{R}_{\ge0}^N\times \mathbb{Z}_{\ge0}^N \\
 &\times\mathbb{Z}_{\ge0}^N\times\mathbb{Z}_{\ge0}^N\times\mathbb{Z}_{\ge0}^N\times \{-1,1\}^N\times \mathbb{R}_{\ge0}^N,\\
 n_{\theta}:=&\sum_{i=1}^{N}n_{\theta,i},
\end{split}\]
will be formulated as a hybrid system model using the formalism in \cite{Dolk Borgers and Heemels,Goebelbook}. Note that the state space $\mathcal{X}$ considers the special structure of $\theta$ defined in (\ref{def of theta}). Meanwhile, define the subspace of $\mathcal{X}$ corresponding to $[\bar{k}_i,l_i,\theta_i^{\rm{T}}]^{\rm{T}}$ as $\tilde{\mathcal{X}}_i$.

Then, the closed-loop system can be described as
\begin{equation}\label{closed loop system}
\begin{cases}
    \dot{q}=F(q,w), &q\in\mathcal{C},\\
  q^+\in G(q), &q\in\mathcal{D},
\end{cases}
\end{equation}
with the flow set $\mathcal{C}$ and jump set $\mathcal{D}$:
\begin{equation}
  \label{flow and jump set}
  \begin{split}
    \mathcal{C}:=&\{q\in\mathcal{X}|\hat{\tau}_i\in[0,T_M^i],i\in\bar{N}\},\\
    \mathcal{D}:=& \cup_{i=1}^N \mathcal{D}_i, \\
    \mathcal{D}_i:=&\{q\in\mathcal{X}|\hat{\tau}_i\in[T_m^i,T_M^i]\}.
  \end{split}
\end{equation}
The flow dynamics is given by
\begin{equation}
  \label{flow closed loop}
  \begin{split}
    F(q,w)=[&f^{\rm{T}}(q,w),g^{\rm{T}}(q,w),0_{n_{\theta}}^{\rm{T}},0_{n_v}^{\rm{T}},\mathbf{1}_{N}^{\rm{T}},\\
    &0_N^{\rm{T}},0_N^{\rm{T}},0_N^{\rm{T}},0_N^{\rm{T}},f_{\eta}^{\rm{T}}(q)]^{\rm{T}}
  \end{split}
\end{equation}
with
\[\begin{split}
  f(q,w):=[&f_p^{\rm{T}}(x_p,g_c(x_c)+e_u,w),f_c^{\rm{T}}(x_c,g_p(x_p)+e_y)]^{\rm{T}}\\
  :=[&\bar{f}_p^{\rm{T}}(q,w),\bar{f}_c(q)]^{\rm{T}},\\
  g(q,w):=\Big[&-\left(\frac{\partial g_p}{\partial x_p}\bar{f}_p(q,w)\right)^{\rm{T}}, -\left(\frac{\partial g_c}{\partial x_c}\bar{f}_c(q)\right)^{\rm{T}}\Big]^{\rm{T}}\\
  f_{\eta}(q):=[&(f_{\eta}^1(q))^{\rm{T}},\dots,(f_{\eta}^N(q))^{\rm{T}}]^{\rm{T}},
\end{split}\]
where $e_y:=\hat{y}-y, e_u:=\hat{u}-u$ and with a little abuse of notation, we express $f_{\eta}^i(o_i,\eta_i)$ as $f_{\eta}^i(q)$ for $i\in\bar{N}$. The set-valued jump map $G(q)$ is defined as $G(q)=\cup_{i=1}^NG_i(q)$ with
\begin{equation}
  \label{jump closed loop}
  G_i(q)=\begin{cases}
    G_i^{1}(q), & m_i=1 \wedge  g_s^i<0 \wedge q\in\mathcal{D}_i,\\
    G_i^{2}(q), & m_i=1 \wedge g_s^i>0 \wedge q\in\mathcal{D}_i,\\
    \{G_i^{1}(q),G_i^{2}(q)\}, &m_i=1 \wedge g_s^i=0\wedge q\in\mathcal{D}_i ,\\
    G_i^{3}(q), & m_i=-1 \wedge q\in\mathcal{D}_i,\\
    \emptyset, & q\notin \mathcal{D}_i,
  \end{cases}
\end{equation}
where the triggering function $g_s^i$ is defined in (\ref{ET general}) and the set-valued jump maps or functions $G_i^n,n=1,2,3,$ are given by
\begin{equation}
  \label{def of Gij1}
  \begin{split}
  G_i^1:=&\left[
    x^{\rm{T}}, e^{\rm{T}}, \left[\begin{smallmatrix}
      \theta_1\\
      \vdots\\
      \theta_{i-1}\\
      \left[\begin{smallmatrix}
        \theta_{i,1}\\
        \vdots\\
        \theta_{i,l}\\
        h_{v_i}-e_i-\sum_{j=1}^{l}\theta_{i,j}\\
        0\\
        \vdots\\
        0
      \end{smallmatrix}\right]\\
      \theta_{i+1}\\
      \vdots\\
      \theta_{N}
    \end{smallmatrix}\right]^{\rm{T}},\left[\begin{smallmatrix}
      \tilde{v}_1\\
      \vdots\\
      \tilde{v}_{i-1}\\
      v_i\\
      \tilde{v}_{i+1}\\
      \vdots\\
      \tilde{v}_{N}
    \end{smallmatrix}\right]^{\rm{T}},\right.\\
    &{\text{\quad}}\left.\hat{\tau}^{\rm{T}}\Gamma_i,\bar{k}^{\rm{T}}+\mathbf{1}_{N}^{\rm{T}}(I-\Gamma_i),\tilde{k}^{\rm{T}},l^{\rm{T}}+\mathbf{1}_{N}^{\rm{T}}(I-\Gamma_i),\right.\\
    &{\text{\quad}}\left.\hat{m}^{\rm{T}}\Gamma_i+\{-\mathbf{1}_{N}^{\rm{T}}(I-\Gamma_i),\xi_{D_i}(l_i)\mathbf{1}_{N}^{\rm{T}}(I-\Gamma_i)\},\right.\\
    &{\text{\quad}}\left.\eta^{\rm{T}}\Gamma_i+g_t^{\rm{T}}(I-\Gamma_i)\right]^{\rm{T}}\\
    \end{split}
\end{equation}
    \begin{equation}
  \label{def of Gij2}
  \begin{split}
    G_i^2:=&\left[x^{\rm{T}},e^{\rm{T}}, \theta^{\rm{T}}, \left[\begin{smallmatrix}
      \tilde{v}_1\\
      \vdots\\
      \tilde{v}_{i-1}\\
      v_i\\
      \tilde{v}_{i+1}\\
      \vdots\\
      \tilde{v}_{N}
    \end{smallmatrix}\right]^{\rm{T}}, \hat{\tau}^{\rm{T}}\Gamma_i,\bar{k}^{\rm{T}}, \tilde{k}^{\rm{T}},l^{\rm{T}},\hat{m}^{\rm{T}}, \right.\\
    &{\text{\quad}}\left. (\Gamma_i\eta+(I-\Gamma_i)g_s)^{\rm{T}}\right]^{\rm{T}}\\
    \end{split}
\end{equation}
    \begin{equation}
  \label{def of Gij3}
  \begin{split}
    G_i^3:=&\left[x^{\rm{T}},\left[\begin{smallmatrix}
      e_1\\
      \vdots\\
      e_{i-1}\\
      e_i+\theta_{i,1}\\
      e_{i+1}\\
      \vdots\\
      e_N
    \end{smallmatrix}\right]^{\rm{T}}, \left[\begin{smallmatrix}
      \theta_1\\
      \vdots\\
      \theta_{i-1}\\
      \left[\begin{smallmatrix}
        \theta_{i,2}\\
        \vdots\\
        \theta_{i,D_i+1}\\
        0
      \end{smallmatrix}\right]\\
      \theta_{i+1}\\
      \vdots\\
      \theta_{N}
    \end{smallmatrix}\right]^{\rm{T}},  \tilde{v}^{\rm{T}},  \hat{\tau}^{\rm{T}}, \bar{k}^{\rm{T}}, \right.\\
    &{\text{\quad}}\left.\tilde{k}^{\rm{T}}+\mathbf{1}_{N}^{\rm{T}}(I-\Gamma_i), l^{\rm{T}}-\mathbf{1}_{N}^{\rm{T}}(I-\Gamma_i),\right.\\
    &{\text{\quad}}\left.\hat{m}^{\rm{T}}\Gamma_i+\{-\xi_{1}(l_i)\mathbf{1}_{N}^{\rm{T}}(I-\Gamma_i),\mathbf{1}_{N}^{\rm{T}}(I-\Gamma_i)\},\eta^{\rm{T}}
   \right]^{\rm{T}},
\end{split}
\end{equation}
where $h_{v_i}$ denotes $h_{v_i}(\bar{k}_i,e_i+\sum_{j=1}^{l}\theta_{i,j})$ and
\[\begin{split}
  g_t:=[g_t^1,\dots,g_t^N]^{\rm{T}}, {\text{\quad}} g_s:=[g_s^1,\dots,g_s^N]^{\rm{T}}.
\end{split}\]
The matrix $\Gamma_i\in\mathbb{R}^{N\times N}$ is a diagonal matrix with its diagonal elements being 1 except the $i$-th element which is $0$. The scalar function $\xi_m(p)=-1$ if $m=p$; otherwise $\xi_m(p)=1$.

The the set-valued jump maps $G_i^n,n=1,2,3,$ describe how the state $q$ jumps when the communication channel $\mathcal{C}_i$ conducts the action of transmission, sampling, and updating, respectively. Note that the evolution of $\theta$ in different cases follows the properties in Proposition \ref{pps of thetai}. Moreover, referring to \cite{Wangperiodic}, the union form $\{G_i^{1}(q),G_i^{2}(q)\}$ in the case of $g_s^i=0$ is used to ensure the outer semi-continuity of $G(q)$ and the resultant nominal well-posedness for the hybrid system, see \cite{Goebelbook} for more details.

\begin{rmk}\label{rmk of m and l}
  In the small-delay case, from the fact that transmissions and updates occur in turn,  \cite{Heemels2010} only introduced $l_i(t)$, which could replace the role of $\hat{m}_i(t)$ and must switch between $0$ and $1$ in turn. However, in the more general large-delay case, the value of $\hat{m}_i(t)$ cannot be decided by the current action. Thus, in (\ref{def of Gij1}) and (\ref{def of Gij3}), the value of $\hat{m}_i(t)$ after a transmission or updating action is given by a set. Moreover, to avoid $q$ jumping out of the definition set $\mathcal{X}$, one has that $\hat{m}_i=1$ when $l_i=0$ and $\hat{m}_i=-1$ when $l_i=D_i+1$ because, respectively, an update must occur after its corresponding transmission, and one transmitted signal must arrive at the destination before experiencing $D_i+1$ subsequent transmissions from Assumption \ref{asm of large delay}. The analysis above leads to the introduction of the function $\xi_{(\cdot)}(\cdot)$ in (\ref{def of Gij1}) and (\ref{def of Gij3}).
\end{rmk}

\subsection{Study Objective}
The interest of this paper is to design the decentralized PETC in (\ref{ET general}), namely, constant $T_M^i$ and functions $f_\eta^i,g_t^i,g_s^i$, for all communication channels $i\in\bar{N}$ to ensure the input-to-state stability, with the definition given as follows \cite{Wangperiodic}.

\begin{dfn}\label{definition of ISS}
  For the closed-loop system in (\ref{closed loop system}-\ref{def of Gij3}), the set
  \[\mathcal{S}:=\{q\in\mathcal{X}|x=0,e=0,\eta=0\} \]
   is input-to-state stable with respect to $w$ if there exist $\beta\in\mathcal{KL}$ and $\psi\in\mathcal{K}_{\infty}$ such that any solution pair $(q,w)$ satisfies\footnote{The definition of $\left\|w\right\|_{\infty}$ in a hybrid system sense can be found in \cite{CaiLetter}.}
  \begin{equation}\label{inequality of ISS}
    \left\|q(t,\bar{j})\right\|_{\mathcal{S}}\le \beta(\left\|q(0,0)\right\|_{\mathcal{S}},t+\bar{j})+\psi(\left\|w\right\|_{\infty}),
  \end{equation}
  for all $(t,\bar{j})\in {\text{dom }} q$, where the initial state $q(0,0)$ satisfies $\tilde{v}_i(0,0)=g_{v_i}(x(0,0))$ for all $i\in\bar{N}$.
\end{dfn}

Due to the trivial assumption of $s_0^i=0$ for all $i\in\bar{N}$, we consider the restriction, $\tilde{v}_i(0,0)=g_{v_i}(x(0,0))$, of the initial state in Definition \ref{definition of ISS}.

\section{Main results}
In this section, the main results of this paper are provided. First, according to some assumptions on the storage functions of the system state and updating errors, sufficient conditions on the network setup and dynamic PETC are proved to ensure the input-to-state stability. Then, the construction of these storage functions are given based on some generally acceptable conditions on storage functions \cite{Heemels2010} from delay-free cases. Finally, the implementation of the proposed dynamic PETC is discussed under different capacities of equipment in communication channels.

\subsection{Stability Analysis}
The stability analysis is given by starting from the following assumptions.
\begin{asm}\label{asm of we}
  For each $i\in\bar{N}$, there exist a function $\tilde{W}_i: \tilde{\mathcal{X}}_i\times \mathbb{R}^{n_{v,i}} \to\mathbb{R}_{\ge0}$ with $\tilde{W}_i(\bar{k}_i,l_i,\cdot,\cdot)$ locally Lipschitz for all fixed $\bar{k}_i\in\mathbb{Z}_{\ge0}$ and $l_i\in\{0,\dots,D_i+1\}$, $\mathcal{K}_{\infty}$--functions $\underline{\beta}_{\tilde{W}_i}$ and $\bar{\beta}_{\tilde{W}_i}$, continuous functions $H_{l_i,i}:\mathbb{R}^{n_x}\times\mathbb{R}^{n_v}\times \mathbb{R}^{n_w}\to \mathbb{R}_{\ge0}$, positive constants $L_{l_i,i}$ for $l_i\in\{0,\dots,D_i+1\}$, and a scalar $0<\tilde{\lambda}_i<1$ such that, for all $[\bar{k}_i,l_i,\theta_i^{\rm{T}}]^{\rm{T}}\in \tilde{\mathcal{X}}_i$, the following statements hold:
\begin{equation}\label{positive of wi}
  \underline{\beta}_{\tilde{W}_i}(\left\|[e_i^{\rm{T}},\theta_i^{\rm{T}}]^{\rm{T}}\right\|)\le \tilde{W}_i(\bar{k}_i,l_i,\theta_i, e_i) \le \bar{\beta}_{\tilde{W}_i}(\left\|[e_i^{\rm{T}},\theta_i^{\rm{T}}]^{\rm{T}}\right\|);
\end{equation}
\begin{subequations}\label{jump of wi}
\begin{equation}
\begin{split}
  \tilde{W}_i&(\bar{k}_i+1,l_i+1,[\theta_{i,1}^{\rm{T}},\dots,\theta_{i,l_i}^{\rm{T}},h_{v_i}^{\rm{T}}-e_i^{\rm{T}}-\sum_{j=1_i}^{l}\theta_{i,j}^{\rm{T}},0^{\rm{T}}]^{\rm{T}}, e_i)\\
  &\le \tilde{\lambda}_i \tilde{W}_i(\bar{k}_i,l_i,\theta_i, e_i),
\end{split}
\end{equation}
\begin{equation}
\begin{split}
  \tilde{W}_i&(\bar{k}_i,l_i-1,[\theta_{i,2}^{\rm{T}},\dots,\theta_{i,D_i+1}^{\rm{T}},0^{\rm{T}}]^{\rm{T}}, e_i+\theta_{i,1})\\
  &\le \tilde{W}_i(\bar{k}_i,l_i,\theta_i, e_i);
\end{split}
\end{equation}
\end{subequations}
for all $e_i\in \mathbb{R}^{n_{v,i}}$; and
\begin{equation}
  \label{flow of wi}
  \left\langle\frac{\partial \tilde{W}_i}{\partial e_i}, -f_{v_i}(q,w)\right\rangle\le L_{l_i,i}\tilde{W}_i+H_{l_i,i}(x,e,w),
\end{equation}
for almost all $e_i\in\mathbb{R}^{n_{v,i}}$, where $f_{v_i}:=\frac{\partial g_{v_i}}{\partial x}f(q,w)$ is the $i$-th block of $-g(q,w)$ in (\ref{flow closed loop}) corresponding to $e_i$, and $h_{v_i}$ denotes $h_{v_i}(\bar{k}_i,e_i+\sum_{j=1}^{l}\theta_{i,j})$.
\end{asm}

{Assumption \ref{asm of we} supposes a storage function $\tilde{W}_i$ on updating error $e_i$ for $i\in\bar{N}$, the jump behavior of which is characterized by (\ref{jump of wi}). That is, by (\ref{jump of wi}a), $\tilde{W}_i$ decays in the rate of $\tilde{\lambda}_i$ when experiencing a transmission; and by (\ref{jump of wi}b), $\tilde{W}_i$ never increases whenever the updating action occurs.  Meanwhile, from (\ref{flow of wi}), $\tilde{W}_i$ has an exponential growth rate $L_{l_i,i}$ in flow.} In the next subsection, we will show how to construct $\tilde{W}_i$ from some general acceptable conditions on the protocol $h_{v_i}$ used for the delay-free NCSs. Recall that, for $[\bar{k}_i,l_i,e_i^{\rm{T}},\theta_i^{\rm{T}}]^{\rm{T}}$, its subspace $\tilde{\mathcal{X}}_i$ depends on the special structure of $\theta_i$ in (\ref{def of theta}), that is, $\theta_{i,j}=0$ with $j=l_i+1,\dots,D_i+1$ for all $i\in\bar{N}$.

\begin{asm}\label{asm of V}
  There exist a locally Lipschitz function $\tilde{V}:\mathbb{R}^{n_x}\to\mathbb{R}_{\ge0}$, locally Lipschitz functions $\tilde{\delta}_i,\hat{\delta}_i:\mathbb{R}^{n_{v,i}}\to\mathbb{R}_{\ge0}$ satisfying $\tilde{\delta}_i(0)=0$, $\mathcal{K}_{\infty}$--functions $\underline{\beta}_{\tilde{V}}, \bar{\beta}_{\tilde{V}}, \underline{\beta}_{\tilde{\delta}_i}, \bar{\beta}_{\tilde{\delta}_i} \alpha_{\tilde{V}},\alpha_w,\sigma_{l_i,i}$, continuous functions $\tilde{J}_i:\mathbb{R}^{n_x}\times\mathbb{R}^{n_v}\times \mathbb{R}^{n_w}\to \mathbb{R}_{\ge0}$, and scalars $\gamma_{l_i,i}>0,\epsilon_{l_i,i}\ge0$, where $i\in\bar{N}$ and $l_i\in\{0,\dots,D_i+1\}$, such that, for all $x\in\mathbb{R}^{n_x}, v_i\in\mathbb{R}^{n_{v_i}}$ and $i\in\bar{N}$,
  \begin{equation}\label{positive of V}
  \begin{split}
    \underline{\beta}_{\tilde{V}}(\left\|x\right\|)\le \tilde{V}(x) \le \bar{\beta}_{\tilde{V}}(\left\|x\right\|),\\
    \underline{\beta}_{\tilde{\delta}_i}(\left\|v_i\right\|)\le \tilde{\delta}_i(v_i) \le \bar{\beta}_{\tilde{\delta}_i}(\left\|v_i\right\|);
  \end{split}
\end{equation}
for almost all $x\in\mathbb{R}^{n_x}$, all $[v^{\rm{T}},e^{\rm{T}},w^{\rm{T}}]^{\rm{T}}\in\mathbb{R}^{2n_v+n_w}$,
\begin{equation}
  \label{flow of V}
  \begin{split}
     \left\langle \nabla \tilde{V}(x), f(q,w)\right\rangle\le& -\alpha_{\tilde{V}}(\left\|x\right\|)+\sum_{i=1}^N\Big(\gamma_{l_i,i}^2\tilde{W}_i^2-\tilde{\delta}_i(v_i)\\
     &-\sigma_{l_i,i}(\tilde{W}_i)-(1+\epsilon_{l_i,i})H_{l_i,i}^2(x,e,w)\\
     &-\hat{\delta}_i(\hat{v}_i)-(1+\epsilon_{l_i,i})\tilde{J}_i(x,e,w)\Big)\\
     &+\alpha_w(\left\|w\right\|);
  \end{split}
\end{equation}
and for almost all $v_i\in\mathbb{R}^{n_{v,i}}$, all $[x^{\rm{T}},e^{\rm{T}},w^{\rm{T}}]^{\rm{T}}\in\mathbb{R}^{n_x+n_v+n_w}$, and all $i\in\bar{N}$,
\begin{equation}
  \label{flow of deltai}
  \begin{split}
     \left\langle \nabla \tilde{\delta}_i(v_i), f_{v_i}(q,w)\right\rangle\le& H_{l_i,i}^2(x,e,w)+\tilde{J}_i(x,e,w),
  \end{split}
\end{equation}
where $\hat{v}_i=v_i+e_i$ is the updated signal in the destination node, and the arguments of $\tilde{W}_i$ are omitted for simplicity.
\end{asm}

{Some analysis on these assumptions are given as follows. }
\begin{itemize}
  \item {The relationship in (\ref{positive of V}) and (\ref{flow of V})}  means that the closed-loop system in (\ref{closed loop system}--\ref{def of Gij3}) is input-to-state stable with respect to $w$ if $e(t)=0$ for all $t\ge0$, since the assumptions are reduced to $\left\langle \nabla \tilde{V}(x), f(q,w)\right\rangle\le -\alpha_{\tilde{V}}(\left\|x\right\|)+\alpha_w(\left\|w\right\|)$. This comes from the fact that the controller in (\ref{controller}) is designed in an emulation-based manner. However, they do not not imply any stability of the closed-loop system since $e$ is an internal state in $q$.
  \item {Motivated by \cite{Wangperiodic}, the equation in (\ref{flow of deltai}) related to $v_i$ is used to characterize some local properties of channel $\mathcal{C}_i$, which ensure the decentralized nature of ETC in each channel.}  The function $\tilde{J}_i(x,e,w)$ stands for the redundant terms in the derivative of $\tilde{\delta}_i(v_i)$ except $H_{l_i,i}^2(x,e,w)$.
  \item  In Assumption~\ref{asm of V}, the terms $\epsilon_{l_i,i}$ and $\hat{\delta}_i$, for $i\in\bar{N}$ and $l_i\in\{0,\dots,D_i+1\}$, are used to facilitate the design of event-triggering conditions (which will be illustrated later). Note that these two terms are not restrictive since they can be selected as zero functions/scalars, although this could lead to more frequent transmissions. {Especially, we consider $\hat{\delta}_i(\hat{v}_i)$ in the flow of $\tilde{V}_i$ because the other signals in (\ref{flow of V}) require continuous reading during sampling instants, which is forbidden for the ET in Fig.~\ref{configuration} according to our setup.}
\end{itemize}

To cope with the network-induced updating error $e$, we introduce a group of variables $\phi_{l_i,i}\in\mathbb{R}_{>0}$ for $i\in\bar{N}$ and $l_i\in\{0,\dots,D_i+1\}$, whose evolution is given by
\begin{equation}\label{def of phi}
  \dot{\phi}_{l_i,i}=-2L_{l_i,i}\phi_{l_i,i}-\gamma_{l_i,i}(\phi_{l_i,i}^2+1),
\end{equation}
where the constants $L_{l_i,i}$ and $\gamma_{l_i,i}$ are from Assumptions \ref{asm of we}--\ref{asm of V}, and the corresponding boundary conditions will be given in the next theorem for characterizing the network setup. For all $i\in\bar{N}$, define the constants
\begin{equation}\label{bound of gamma}
  \underline{\gamma}_i:=\min_{l_i\in\{0,\dots,D_i+1\}}\{\gamma_{l_i,i}\}, {\text{\quad}} \bar{\gamma}_i:=\max_{l_i\in\{0,\dots,D_i+1\}}\{\gamma_{l_i,i}\},
\end{equation}
and the variable $\varpi_i\in\mathbb{R}_{>0}$, which evolves according to
\begin{equation}\label{def of varpi}
  \dot{\varpi}_i(t)=-\frac{1-\pi_i}{T_M^i}, {\text{\quad}} \varpi_i(0)=1,
\end{equation}
with some free parameter $\pi_i\in(0,1)$, which results in $\varpi_i(\hat{\tau}_i)>0$ for all $i\in\bar{N}$.

\begin{thm}
  \label{thm of stability}
  For the closed-loop system in (\ref{closed loop system}--\ref{def of Gij3}) under Assumptions \ref{asm of large delay}--\ref{asm of V}, the set $\mathcal{S}$ in Definition \ref{definition of ISS} is input-to-state stable with respect to $w$ if the network setup is given as follows:
  \begin{enumerate}
    \item for each communication channel $\mathcal{C}_i,i\in\bar{N}$, the upper bound $T_M^i$ and initial conditions, $\phi_{l_i,i}(0)>0$, satisfy
    \begin{equation}\label{condition TMI transmission}
\gamma_{l_i,i}\phi_{l_i,i}(T_M^i)\ge \tilde{\lambda}_i^2\gamma_{l_i+1,i}\phi_{l_i+1,i}(0),
\end{equation}
 for all $l_i\in\{0,\dots, D_i\}$; and
     \begin{equation}\label{condition TMI update}
\gamma_{l_i-1,i}\phi_{l_i-1,i}(\hat{\tau}_i)\le \gamma_{l_i,i}\phi_{l_i,i}(\hat{\tau}_i),
\end{equation}
for all $\hat{\tau}_i\in[0,T_M^i]$ and $l_i\in\{1,\dots, D_i+1\}$;
    \item for each communication channel $\mathcal{C}_i,i\in\bar{N}$, the functions in the dynamics of $\eta_i$ in (\ref{dynamics of etai}) satisfy
   \begin{subequations}\label{conditions on dynamics of eta}
\begin{equation}
\begin{split}
 f_{\eta}^i(q)\le& -a_i\eta_i+\hat{\delta}_i(\hat{v}_i)+(1-\varepsilon_i)\tilde{\rho}_i\tilde{\delta}_i(\tilde{v}_i),
\end{split}
\end{equation}
\begin{equation}
\begin{split}
g_t^i(q)\le&\eta_i-\bar{\rho}_i\tilde{\delta}_i(v_i)+\bar{\rho}_i\max\{\tilde{\delta}_i(v_i),\varpi_i(\hat{\tau}_i)\tilde{\delta}_i(\tilde{v}_i)\}\\
  &-\max\{\gamma_{l_i+1,i}\phi_{l_i+1,i}(0)\tilde{\lambda}_i^2\tilde{W}_i^2,\hat{\rho}_i\tilde{\delta}_i(v_i)\}\\
  &+\max\{\gamma_{l_i,i}\phi_{l_i,i}(\hat{\tau}_i)\tilde{W}_i^2,\hat{\rho}_i\tilde{\delta}_i(v_i)\},
\end{split}
\end{equation}
\begin{equation}
\begin{split}
  g_s^i(q)\le&\eta_i-\bar{\rho}_i\tilde{\delta}_i(v_i)+\bar{\rho}_i\max\{\tilde{\delta}_i(v_i),\varpi_i(\hat{\tau}_i)\tilde{\delta}_i(\tilde{v}_i)\}\\
  &-\max\{\gamma_{l_i,i}\phi_{l_i,i}(0)\tilde{W}_i^2,\hat{\rho}_i\tilde{\delta}_i(v_i)\}\\
  &+\max\{\gamma_{l_i,i}\phi_{l_i,i}(\hat{\tau}_i)\tilde{W}_i^2,\hat{\rho}_i\tilde{\delta}_i(v_i)\},
\end{split}
\end{equation}
\end{subequations}
for all $l_i\in\{0,\dots,D_i+1\}$, where $a_i>0$ and $\varepsilon_i\in(0,1)$ are user-specified free parameters and the constants $\bar{\rho}_i,\tilde{\rho}_i,\hat{\rho}_i\ge0$ satisfy
\begin{equation}
  \nonumber
  \begin{split}
    \bar{\rho}_i\le&\min_{l_i\in\{0,\dots,D_i+1\}}\{\epsilon_{l_i,i}\},\\
    \hat{\rho}_i\le&\frac{1}{2}\min\{1,\underline{\phi}_i/\bar{\gamma}_i\},\\
    \tilde{\rho}_i:=&\min\left\{\frac{\bar{\rho}_i(1-\pi_i)}{T_M^i},\frac{\pi_i}{2}\right\},
  \end{split}
\end{equation}
with $\epsilon_{l_i,i}$ defined in Assumption \ref{asm of V} and $\underline{\phi}_i,\bar{\phi}_i,\underline{\varpi}_i>0$ ensuring
\begin{equation}\label{bound of phi}
  \begin{split}
  \underline{\phi}_i\le& \phi_{l_i,i}(\hat{\tau}_i) \le\bar{\phi}_i,{\text{\quad}}\underline{\varpi}_i \le \varpi_i(\hat{\tau}_i)\le 1,
\end{split}
\end{equation}
for all $\hat{\tau}_i\in[0,T_M^i]$, $l_i\in\{0,\dots,D_i+1\}$ and $i\in\bar{N}$.
  \end{enumerate}
\end{thm}

\begin{IEEEproof}
  See Appendix B.
\end{IEEEproof}

In the flow dynamics of $\eta_i$, the function $f_{\eta}^i$ only depends on, besides $\eta_i$, the sampled and scheduled signals, $\tilde{v}_i$ and $\hat{v}_i$. This agrees with the sampled-data structure in Fig.~\ref{configuration}.
\begin{rmk}
  \label{rmk of free parameters}
{Combining (\ref{ET general}) and (\ref{conditions on dynamics of eta}c) implies that the increase of $\eta_i$ helps in reducing the number of events; while from (\ref{flow of U}) in Appendix B one has that the convergence rate of closed-loop systems is affected by $a_i$ and $\varepsilon_i$. Thus, these two free parameters introduce a tradeoff between transmission and stability performance: smaller $a_i$ and $\varepsilon_i$ generate less events but slow down the convergence.}
\end{rmk}
\begin{rmk}\label{rmk of static ET}
  The static event-triggering condition, which is independent of $\eta_i$, can be described as
\begin{equation}
  \label{general ET static}
  t_{k+1}^i=\min\{t>t_k^i|t\in\{s_j^i\}_{j\in\mathbb{Z}_{\ge0}}, g_s^i(o_i(t),0)<0\},
\end{equation}
Note that (\ref{general ET static}) can be analyzed from Theorem \ref{thm of stability} by designing $f_{\eta}^i,g_t^i$ and $g_s^i$ to ensure $\eta_i(t)=0$ for all $t\ge0$ and $i\in \bar{N}$. In \cite{Wangperiodic}, an event-triggering condition was designed by replacing $g_s^i(o_i,0)<0$ in (\ref{general ET static}) as $\gamma_{l_i,i}\phi_{l_i,i}(0)\tilde{W}_i^2>\hat{\rho}_i\tilde{\delta}_i(v_i)$, which is more conservative because $\gamma_{l_i,i}\phi_{l_i,i}(0)\tilde{W}_i^2\le\hat{\rho}_i\tilde{\delta}_i(v_i)$ can render the right-hand side of (\ref{conditions on dynamics of eta}c) to be nonnegative.
\end{rmk}

\begin{rmk}\label{rmk of time triggered}
{Considering the event-triggering condition in (\ref{general ET static}) with $\bar{\rho}_i=0$ and $\hat{\rho}_i=0$, one has that $g_s^i(o_i(t),0)$ is trivially negative. Hence, only Item 1) in Theorem \ref{thm of stability} works and it can cover the case of time-triggering control, where each sampling instant corresponds to a transmission. Especially, in the small-delay case, Item 1) in Theorem \ref{thm of stability} is almost the same as (25) in \cite{Heemels2010}. One slight difference is that the maximum length of delays could be smaller than the sampling period in \cite{Heemels2010} while by setting $D_i=0$, we suppose that the worst delay is equal to the sampling period. Therefore, Theorem \ref{thm of stability} generalizes the framework in [3] to the large-delay case.}
\end{rmk}

\begin{rmk}\label{rmk of extra assumption}
  Without the extra terms $\hat{\delta}_i$ and $\epsilon_{l_i,i}$ in Assumption \ref{asm of V} and the introduction of $\varpi_i$, the conditions in (\ref{conditions on dynamics of eta}) become $f_{\eta}^i(q)\le -a_i\eta_i$, and
  \[\begin{split}
    g_s^i(q)\le&\eta_i-\max\{\gamma_{l_i,i}\phi_{l_i,i}(0)\tilde{W}_i^2,\hat{\rho}_i\tilde{\delta}_i(v_i)\}\\
  &+\max\{\gamma_{l_i,i}\phi_{l_i,i}(\hat{\tau}_i)\tilde{W}_i^2,\hat{\rho}_i\tilde{\delta}_i(v_i)\}\\
  \le& \eta_i,
  \end{split}\]
  which imply that $\eta_i(t)$ is unable to increase during two consecutive transmission instants. In this case, the dynamic event-triggering condition will approximately reduce to a static one and generate more events.
\end{rmk}

\subsection{Construction of  Functions}
The functions, $\tilde{W}_i$ and $\tilde{V}$, $i\in\bar{N}$, play an important role in Theorem \ref{thm of stability}, by characterizing the impacts of delayed signals on the system behavior. Thus, we will show how to construct them from the conditions used for the delay-free cases. {Referring to \cite{Wangperiodic}--\cite{Dolk Borgers and Heemels}, we introduce the following basic assumptions.}

\begin{asm}\label{asm of protocol}
For each $i\in\bar{N}$, there exist a function ${W}_i:\mathbb{Z}_{\ge0}\times \mathbb{R}^{n_{v,i}} \to\mathbb{R}_{\ge0}$ with $W_i(\bar{k}_i,\cdot)$ continuous for all fixed $\bar{k}_i\in\mathbb{Z}_{\ge0}$, and constants $\bar{\beta}_{{W}_i}\ge\underline{\beta}_{{W}_i}>0,{M_{p,i}}\ge0,\lambda_{W,i}\ge 1,0\le\lambda_i<1$ such that, for all $\bar{k}_i\in\mathbb{Z}_{\ge0}$, the following statements hold:
\begin{enumerate}
  \item $\underline{\beta}_{W_i}\left\|e_i\right\|\le W_i(\bar{k}_i, e_i) \le \bar{\beta}_{W_i}\left\|e_i\right\|$ for all $e_i\in\mathbb{R}^{n_{v,i}}$;
  \item $W_i(\bar{k}_i+1,h_{v_i}(\bar{k}_i,e_i))\le\lambda_i W_i(\bar{k}_i,e_i)$ for all $e_i\in\mathbb{R}^{n_{v,i}}$;
  \item $W_i(\bar{k}_i+1,e_i)\le \lambda_{W,i} W_i(\bar{k}_i,e_i)$ for all $e_i\in\mathbb{R}^{n_{v,i}}$;
  \item $\left\|\frac{\partial W_i}{\partial e_i}(\bar{k}_i,e_i)\right\|\le M_{p,i}$  for almost all $e_i\in\mathbb{R}^{n_{v,i}}$,
\end{enumerate}
where $h_{v_i}$ is the scheduling protocol defined in (\ref{scheduling protocol}).
\end{asm}

\begin{asm}\label{asm of fvi}
  There exist a continuous functions $H_i:\mathbb{R}^{n_x}\times\mathbb{R}^{n_v}\times \mathbb{R}^{n_w}\to \mathbb{R}_{\ge0}$ and a constant $M_{e,i}\ge0$ such that
  \[\left\|f_{v_i}\right\|\le M_{p,i}^{-1}\left(H_i(x,e,w)+M_{e,i}\left\|e_i\right\|\right).\]
\end{asm}

\begin{asm}\label{asm of V orignal}
  There exist a locally Lipschitz function $V:\mathbb{R}^{n_x}\to\mathbb{R}_{\ge0}$, locally Lipschitz functions $\delta_i,\bar{\delta}_i:\mathbb{R}^{n_{v,i}}\to\mathbb{R}_{\ge0}$ satisfying $\delta_i(0)=0$, $\mathcal{K}_{\infty}$--functions $\underline{\beta}_{\tilde{V}}, \bar{\beta}_{\tilde{V}}, \underline{\beta}_{\tilde{\delta}_i}, \bar{\beta}_{\tilde{\delta}_i}, \alpha_{\tilde{V}},\alpha_w,\sigma_{l_i,i}$, continuous functions $J_i:\mathbb{R}^{n_x}\times\mathbb{R}^{n_v}\times \mathbb{R}^{n_w}\to \mathbb{R}_{\ge0}$, and scalars $\gamma_{i}>0,\epsilon_{i}\ge0$, for all $i\in\bar{N}$, such that the following statements hold:
  \begin{enumerate}
    \item for all $x\in\mathbb{R}^{n_x}, v_i\in\mathbb{R}^{n_{v_i}}$, $\underline{\beta}_{V}(\left\|x\right\|)\le V(x) \le \bar{\beta}_{V}(\left\|x\right\|)$ and $\underline{\beta}_{\delta_i}(\left\|v_i\right\|)\le \delta_i(v_i) \le \bar{\beta}_{\delta_i}(\left\|v_i\right\|);$
    \item for almost all $x\in\mathbb{R}^{n_x}$ and all $[v^{\rm{T}},e^{\rm{T}},w^{\rm{T}}]^{\rm{T}}\in\mathbb{R}^{2n_v+n_w}$,
    \begin{equation}
  \nonumber
  \begin{split}
     \left\langle \nabla V(x), f(q,w)\right\rangle\le& -\alpha_{V}(\left\|x\right\|)+\sum_{i=1}^N\Big((\gamma_{i}^2-\bar{\varepsilon}_i)W_i^2\\
     &-\delta_i(v_i)-(1+\epsilon_{i})H_{i}^2(x,e,w)\\
     &-\bar{\delta}_i(\hat{v}_i)-(1+\epsilon_{i}){J}_i(x,e,w)\Big)\\
     &+\bar{\alpha}_w(\left\|w\right\|),
  \end{split}
\end{equation}
    \item for almost all $v_i\in\mathbb{R}^{n_{v,i}}$, all $[x^{\rm{T}},e^{\rm{T}},w^{\rm{T}}]^{\rm{T}}\in\mathbb{R}^{n_x+n_v+n_w}$, and all $i\in\bar{N}$,
    \begin{equation}
    \nonumber
     \left\langle \nabla \delta_i(v_i), f_{v_i}(q,w)\right\rangle\le H_{i}^2(x,e,w)+J_i(x,e,w),
  \end{equation}
  \end{enumerate}
  where the arguments of $W_i$ are omitted for simplicity, and the (sufficiently small) constant $\bar{\epsilon}_i>0$ satisfies $\bar{\varepsilon}_i\in(0,\gamma_{i}^2)$ for all $i\in\bar{N}$.
\end{asm}

\begin{rmk}
  \label{rmk of Assumptions 6}
  {Assumptions \ref{asm of protocol} and \ref{asm of fvi} are the same as (38--42) in \cite{Heemels2010} and (48--52) in \cite{Dolk Borgers and Heemels}. For Assumption \ref{asm of V orignal}, it is reduced to (43-44) in \cite{Heemels2010} by tailoring some terms and conditions for PETC, such as those related with $\delta_i(v_i), \bar{\delta}_i(\hat{v}_i)$ and $J_i$. Meanwhile, if removing $\bar{\delta}_i(\hat{v}_i)$ and $\epsilon_i$ that facilitate the design of dynamic PETC, Assumption \ref{asm of V orignal} is the same as Assumption 2 in \cite{Wangperiodic}. In summary, the assumptions in delay-free cases are standard and more detailed discussions can be found in \cite{Wangperiodic}--\cite{Dolk Borgers and Heemels}.}
\end{rmk}

Then, based on Assumptions \ref{asm of protocol}--\ref{asm of fvi}, we first show that the following form of $\tilde{W}_i(\bar{k}_i,l_i,\theta_i, e_i)$ satisfies Assumption \ref{asm of we}:
\begin{equation}
  \label{def of tilde wi}
  \begin{split}
      \tilde{W}_i(\bar{k}_i,l_i,\theta_i,e_i)
      :=\max\Big\{&\frac{\tilde{\lambda}_i^{l_i}W_i(\bar{k}_i,e_i)}{\lambda_{W,i}^{l_i}},\frac{\tilde{\lambda}_i^{l_i}W_i(\bar{k}_i,e_i+\theta_{i,1})}{\tilde{\lambda}_i\lambda_{W,i}^{l_i-1}},\\
      &\dots,W_i(\bar{k}_i,e_i+\sum_{j=1}^{l_i}\theta_{i,j})\Big\},
  \end{split}
\end{equation}
{where the design parameter $\tilde{\lambda}_i\in[0,1]$ will be specified later.} Note that the number of terms for the maximum operation in $\tilde{W}_i(\bar{k}_i,l_i,\theta_i,e_i)$ is equal to $l_i+1$, rather than $D_i+1$.

\begin{pps}\label{pps of we}
  Suppose that Assumptions \ref{asm of protocol}--\ref{asm of V orignal} hold for the closed-loop system in (\ref{closed loop system}--\ref{def of Gij3}). For each communication channel $\mathcal{C}_i,i\in\bar{N}$, and any $\tilde{\lambda}_i\in(\lambda_i,1)$, the function $\tilde{W}_i(\bar{k}_i,l_i,\theta_i, e_i)$ in (\ref{def of tilde wi}) satisfies all the conditions in Assumption \ref{asm of we} with some $\mathcal{K}_{\infty}$--functions $\underline{\beta}_{\tilde{W}_i}$,$\bar{\beta}_{\tilde{W}_i}$, where the other parameters and functions are given by
  \[L_{l_i,i}=\frac{\lambda_{W,i}^{l_i}M_{e,i}}{\tilde{\lambda}_i^{l_i}\underline{\beta}_{W_i}}, \text{ and } H_{l_i,i}(x,e,w)=H_i(x,e,w)\]
\end{pps}

\begin{IEEEproof}
  See Appendix B.
\end{IEEEproof}

Proportion \ref{pps of we} suggests that a better sensor scheduling protocol (resulting in smaller $\lambda_i$) \cite{Nesic2004} can provide a wider design range of network setup.

\begin{pps}\label{pps of V}
  Suppose that the conditions in Proposition~\ref{pps of we} and Assumption~\ref{asm of V orignal} hold, then the functions $\tilde{V}=V(x)$ and $\tilde{\delta_i}=\delta_i$ satisfy Assumption \ref{asm of V} with the following functions and parameters
  \begin{equation}\label{notations in pps V}
    \begin{split}
       \underline{\beta}_{\tilde{V}}&=\underline{\beta}_{V},\text{ }\bar{\beta}_{\tilde{V}}=\bar{\beta}_{V},\text{ } \underline{\beta}_{\tilde{\delta}_i}=\underline{\beta}_{\delta_i},\text{ } \bar{\beta}_{\tilde{\delta}_i}=\bar{\beta}_{\delta_i},\\
       \gamma_{l_i,i}&=\frac{\lambda_{W,i}^{l_i}\gamma_{i}}{\tilde{\lambda}_i^{l_i}},\text{ } \sigma_{l_i,i}(r)=\frac{\lambda_{W,i}^{2l_i}\bar{\varepsilon}_i}{\tilde{\lambda}_i^{2l_i}}r^2,\text{ }\tilde{J}_i=J_i,\\
       \hat{\delta}_i&=\bar{\delta}_i,\text{ } \alpha_{\tilde{V}}=\alpha_{V},\text{ } \alpha_w=\bar{\alpha}_w,\text{ } \epsilon_{l_i,i}=\epsilon_i,
    \end{split}
  \end{equation}
  for all $l_i\in\{0,\dots,D_i+1\}$ and $i\in\bar{N}$.
\end{pps}

\begin{IEEEproof}
  See Appendix B.
\end{IEEEproof}

\begin{rmk}\label{rmk of derivative of V}
  In \cite{Wangperiodic}, a similar assumption on the derivative of $V$ was given as
  \begin{equation}
    \label{derivative of V reference}
     \begin{split}
     \left\langle \nabla V_{\rm{R}}(x), f(q,w)\right\rangle\le& -\alpha_{{\rm{R}},V}(\left\|x\right\|)+\sum_{i=1}^N\Big((\gamma_{{\rm{R}},i}^2-\bar{\varepsilon}_{{\rm{R}},i})W_i^2\\
     &-\delta_{{\rm{R}},i}(v_i)-H_{R,i}^2-{J}_{R,i}\Big)\\
     &+\bar{\alpha}_{{\rm{R}},w}(\left\|w\right\|),
  \end{split}
  \end{equation}
  where the notations, by adding the subscript ``R'', denote the corresponding version in \cite{Wangperiodic} with respect to Item 2) in Assumption \ref{asm of V orignal}. A systematic design framework was given in \cite{Wangperiodic} to ensure (\ref{derivative of V reference}) for the systems with globally Lipschitz dynamics. Compared to Assumption \ref{asm of V orignal}, there are two differences: $\epsilon_i$ and $\bar{\delta}_i$. To obtain the term $\epsilon_i\ge0$, one can easily select $V:=(1+\epsilon_i) V_{\rm{R}}, H_i=H_{{\rm{R}},i}$ and $J_i= J_{{\rm{R}},i}$, which would lead to an increase on the gain of $W_i^2$ by $\epsilon_i\gamma_{{\rm{R}},i}^2$. Meanwhile, if the function $\delta_{{\rm{R}},i}$ satisfies
  \begin{equation}\label{inequality of check}
      \begin{split}
    \delta_{{\rm{R}},i}(v_i)=&\delta_{{\rm{R}},i}(\hat{v}_i-e_i)\ge \check{\delta}_{{\rm{R}},i}(\hat{v_i})-\check{\epsilon}_{{\rm{R}},i}W_i^2(\bar{k}_i, e_i),
  \end{split}
  \end{equation}
  with some $\check{\delta}_{{\rm{R}},i}:\mathbb{R}^{n_{v_i}}\to\mathbb{R}_{\ge0}$ and $\check{\epsilon}_{{\rm{R}},i}\ge0$, then we have $\bar{\delta}_i:=\epsilon_i\check{\delta}_{{\rm{R}},i}(v_i)$. Subsequently, the introduction of $\epsilon_i$ and $\bar{\delta}_i$ gives $\gamma_i^2=(1+\epsilon_i)\gamma_{{\rm{R}},i}^2+\epsilon_i\check{\epsilon}_{{\rm{R}},i}$. Thus, from Remark~\ref{rmk of extra assumption} and Proposition~\ref{pps of V}, the analysis above implies that the improvement on event-triggering conditions may require a higher sampling frequency. This agrees with the intuition of event-triggering control, that is, collecting more online information via higher sampling frequencies could yield less conservative design of event-triggering conditions. Some similar analysis has been observed in previous studies on PETC, see, e.g., Remark 5 in~\cite{Xiao2020}.
\end{rmk}

\begin{rmk}\label{rmk of check}
  The inequality in (\ref{inequality of check}) holds simply with zero functions and constants; while in some special case, one can have better choices. For example, consider quadratic forms, i.e., $\delta_{{\rm{R}},i}(v_i)=\left\|v_i\right\|^2$ and $W_i^2(\bar{k}_i, e_i)=\left\|e_i\right\|^2$, which are common for the systems with globally Lipschitz nonlinear dynamics and SD protocols \cite{Heijmans}. Then from Young's inequality \cite{Hardy}, we have (\ref{inequality of check}) as
  \[\left\|v_i\right\|^2=\left\|\hat{v}_i-e_i\right\|^2\ge (1-\check{\varepsilon}_i)\left\|\hat{v}_i\right\|^2-(\frac{1}{\check{\varepsilon}_i}-1)\left\|e_i\right\|^2,\]
  with a free parameter $\check{\varepsilon}_i\in(0,1]$. In this way, a similar analysis on $\epsilon_i$ is feasible for $\check{\varepsilon}_i$ as well; that is, a smaller $\check{\varepsilon}_i$ leads to more frequent sampling but better event-triggering conditions. Note that the selection of $\check{\varepsilon}_i$ cannot be arbitrarily small since $T_M^i$ should be larger than the minimum inter-sampling time $T_m^i$, which is decided by hardware constraints in reality.
\end{rmk}

\begin{rmk}\label{rmk of small delay}
   In the small-delay case, namely, $D_i=0$ and $\theta_i=\theta_{i,1}$, \cite{Heemels2010} and \cite{Dolk Borgers and Heemels} designed a storage function $\tilde{W}_i(\bar{k}_i,l_i,\theta_i, e_i)$ as
  \begin{equation}\label{tilde wi existing}
     \begin{cases}
       \max\{W_i(\bar{k}_i,e_i),W_i(\bar{k}_i,e_i+\theta_i)\}, & l_i=0,\\
       \max\{\frac{\tilde{\lambda}_i}{\lambda_{W,i}}W_i(\bar{k}_i,e_i),W_i(\bar{k}_i,e_i+\theta_i)\}, & l_i=1,\\
     \end{cases}
  \end{equation}
  where $\theta_i$ is updated as $-e_i-h_{v_i}$, instead of $0$ as in (\ref{def of theta}), after an updating instant. For (\ref{tilde wi existing}), the extra term $W_i(\bar{k}_i,e_i+\theta_i)$ in the case of $l_i=0$ is used to ensure (\ref{positive of wi}) for all $e_i\in\mathbb{R}^{n_{v_i}}$ and $\theta_i\in\mathbb{R}^{n_{v_i}}$. However, by limiting the consideration of $\theta_i\in\mathbb{R}^{n_{v_i}}$ satisfying $\theta_{i,j}=0$ with $j=l_i+1,\dots,D_i+1$ from (\ref{def of theta}), the extra term $W_i(\bar{k}_i,e_i+\theta_i)$ when $l_i=0$ can be removed.  Thus, the construction of $\tilde{W}_i$ in (\ref{def of tilde wi}) is simpler (in less terms) and more general (in large delays).
\end{rmk}

\subsection{Implementation of Dynamic PETC}
To implement the conditions in (\ref{conditions on dynamics of eta}) in practice, the equipment in the communication channel $\mathcal{C}_i,i\in\bar{N}$, requires the following capacities: (i) the ET is able to solve differential equalities online; (ii) to decide $\hat{v}_i$ and $\gamma_{l_i}\phi_{l_i}$, the holder in the destination node can send back a delay-free acknowledgement signal at every time it receives new updating signals; and (iii) the ET has fundamental capacity to store relevant parameters and realize algebraic and logic operation. In the best case,  where  all the capacities in Items (i)--(iii) are available, one can design the event-triggering condition in (\ref{ET general}) by directly implementing $f_{\eta}^i, g_t^i$, and $g_s^i$ as the corresponding functions in the right-hand side of (\ref{conditions on dynamics of eta}). Otherwise, if the equipment only has limited capacities, we will give the discussions on the implementation of (\ref{ET general}) in this subsection. Note that the capacities of each independent channel is not necessarily the same.

The computational and storage capacity in (iii) is fundamental and necessary. So, in the following, we only consider the case of lacking Item (i) or (ii).

If the ET in $\mathcal{C}_i,i\in\bar{N}$, is unable to solve differential equalities online (namely, lacking Item (i)), it means that the following relationship or signals cannot be obtained in real-time: $\dot{\eta}_i=f_{\eta}^i(o_i(t),\eta_i(t))$ in (\ref{dynamics of etai}), $\hat{\tau}_i$, $\varpi_i(\hat{\tau}_i)$, and $\phi_{l_i,i}(\hat{\tau}_i)$.

Since both $\varpi_i(\cdot)$ and $\phi_{l_i,i}(\cdot)$ are decreasing, one can use $T_M^i$ to give the lower bounds $\varpi_i(T_M^i)$, and $\phi_{l_i,i}(T_M^i)$. For the differential equality in (\ref{dynamics of etai}), note that the terms $\hat{\delta}_i(\hat{v}_i)+(1-\varepsilon_i)\tilde{\rho}_i\tilde{\delta}_i(\tilde{v}_i)$ in (\ref{conditions on dynamics of eta}a) keep constant between two consecutive sampling instants. Thus, one has
\[\begin{split}
  \eta_i(s_{j+1}^i)\le&e^{-a\tau_j^i}\eta_i(s_j^{i+})+\frac{1}{a_{i}}(1-e^{-a_i\tau_j^i})\hat{\delta}_i(\hat{v}_i)\\
  &+\frac{1}{a_{i}}(1-e^{-a_i\tau_j^i})(1-\varepsilon_i)\tilde{\rho}_i\tilde{\delta}_i(\tilde{v}_i)\\
  \le&e^{-aT_M^i}\eta_i(s_j^{i+})+\frac{1}{a_{i}}(1-e^{-a_iT_m^i})\hat{\delta}_i(\hat{v}_i)\\
  &+\frac{1}{a_{i}}(1-e^{-a_iT_m^i})(1-\varepsilon_i)\tilde{\rho}_i\tilde{\delta}_i(\tilde{v}_i).
\end{split}\]
Hence, a conservative event-triggering condition is given by
\begin{equation}\label{implement case i}
\begin{split}
f_{\eta}^i=&0,\\
g_t^i(q)\le&e^{-aT_M^i}\eta_i+\frac{1}{a_{i}}(1-e^{-a_iT_m^i})\hat{\delta}_i(\hat{v}_i)-\bar{\rho}_i\tilde{\delta}_i(v_i)\\
  &+\frac{1}{a_{i}}(1-e^{-a_iT_m^i})(1-\varepsilon_i)\tilde{\rho}_i\tilde{\delta}_i(\tilde{v}_i)\\
  &+\bar{\rho}_i\max\{\tilde{\delta}_i(v_i),\varpi_i(T_M^i)\tilde{\delta}_i(\tilde{v}_i)\}\\
  &-\max\{\gamma_{l_i+1,i}\phi_{l_i+1,i}(0)\tilde{\lambda}_i^2\tilde{W}_i^2,\hat{\rho}_i\tilde{\delta}_i(v_i)\}\\
  &+\max\{\gamma_{l_i,i}\phi_{l_i,i}(T_M^i)\tilde{W}_i^2,\hat{\rho}_i\tilde{\delta}_i(v_i)\},\\
g_s^i(q)\le&e^{-aT_M^i}\eta_i+\frac{1}{a_{i}}(1-e^{-a_iT_m^i})\hat{\delta}_i(\hat{v}_i)-\bar{\rho}_i\tilde{\delta}_i(v_i)\\
  &+\frac{1}{a_{i}}(1-e^{-a_iT_m^i})(1-\varepsilon_i)\tilde{\rho}_i\tilde{\delta}_i(\tilde{v}_i)\\
  &+\bar{\rho}_i\max\{\tilde{\delta}_i(v_i),\varpi_i(T_M^i)\tilde{\delta}_i(\tilde{v}_i)\}\\
  &-\max\{\gamma_{l_i,i}\phi_{l_i,i}(0)\tilde{W}_i^2,\hat{\rho}_i\tilde{\delta}_i(v_i)\}\\
  &+\max\{\gamma_{l_i,i}\phi_{l_i,i}(T_M^i)\tilde{W}_i^2,\hat{\rho}_i\tilde{\delta}_i(v_i)\},
\end{split}
\end{equation}
which shows that, under a given upper bound $T_M^i$, a constant sampling period is better for generating less events due to the increase of $(1-e^{-a_iT_m^i})$.

If there is no  delay-free acknowledgement mechanism (namely, lacking Item (ii)), one cannot obtain the exact values of $l_i$ and the resultant $\theta_i,\gamma_{l_i,i}\phi_{l_i,i}$, $\tilde{W}_i, \hat{v}_i$ and $e_i$.  Hence we design the following left-continuous variable to estimate the upper bound of $l_i(s_{j+1}^i)$. If $D_i=0$, let $\hat{l}_i(t)=0$ for all $t\in\mathbb{R}_{\ge0}$; and if $D_i\in\mathbb{Z}_{\ge1}$, let
\begin{equation}
  \label{def of hat li}
  \begin{split}
  \hat{l}_i(t):=& \hat{l}_i(s_j^{i+})\in \{0,\dots,D_i\},t\in(s_j^i,s_{j+1}^i],\\
  \hat{l}_i(s_j^{i+}):=&\left|\left\{t_k^i\le s_{j}^i|t_k^i=s_m^i,m=j,\dots,j-D_i+1\right\}\right|,
  \end{split}
\end{equation}
which can be obtain based on the information only in the transmitter node.

\begin{pps}\label{pps of upper bound of l}
  $l_i(s_j^i)\le \hat{l}_i(s_j^i)$ holds for any communication $i\in\bar{N}$ and sampling instant $s_j^i,j\in\mathbb{R}_{\ge0}$.
\end{pps}

\begin{IEEEproof}
  See Appendix B.
\end{IEEEproof}

Proposition \ref{pps of upper bound of l} implies that in the small-delay case, we always have that $l_i(s_j^i)=0$ is available, without the introduction of acknowledgement mechanisms.

With the upper bound $\hat{l}_i$ in (\ref{def of hat li}), we can provide the worst-case estimate of $\theta_i,\gamma_{l_i,i}\phi_{l_i,i}$, $\tilde{W}_i$ and $\hat{v}_i$. First consider the following block matrix $\Theta_i=\left[\Theta_{i}^{m,n}\right]_{m,n\in\{1,\dots,D_i+1\}}$ that contains the information of $\theta_i$ with
\[ \Theta_{i}^{m,n}(t)=\bar{v}_i\left(t_{\bar{k}_i(t)-n+m-1}^{i+}\right)-\bar{v}_i\left(t_{\bar{k}_i(t)-n+m-2}^{i+}\right),\]
for $n\in\{2,\dots,D_i+1\}$ and $ m\in\{1,\dots,n\}$; otherwise, $\Theta_{i}^{m,n}(t)=0$. For $e_i(t)$, define the block matrix $E_i=\left[E_i^{m,n}\right]_{m=1,n\in\{1,\dots,D_i+1\}}$ that satisfies
\[ E_{i}^{1,n}(t)=\bar{v}_i\left(t_{\bar{k}_i(t)-n}^{i+}\right)-v_{i},\]
for $n\in\{1,\dots,D_i+1\}$. Thus, if $l_i(s_j^i)=n-1\in\{0,\dots,D_i\}$, we have $e(s_j^i)=E_i^{1,n}(s_j^i)$ and $\theta_i(s_j^i)=\Theta_i^n(s_j^i)$ where $\Theta_i^n$ stands for the $n$-th block column of $\Theta_i$.

Subsequently, for each $n\in\{0,\dots,D_i\}$, define
\[\begin{split}
  \check{v}_{i,n}(t)=&\bar{v}_i\left(t_{\bar{k}_i(t)-n-1}^{i+}\right),\\
  \mathcal{W}_{i,n}^t(\hat{\tau}_i)=&-\max\{\gamma_{n+1,i}\phi_{n+1,i}(0)\tilde{\lambda}_i^2\check{W}_i^2(n),\hat{\rho}_i\tilde{\delta}_i(v_i)\}\\
  &+\max\{\gamma_{n,i}\phi_{n,i}(\hat{\tau}_i)\check{W}_i^2(n),\hat{\rho}_i\tilde{\delta}_i(v_i)\},\\
  \mathcal{W}_{i,n}^s(\hat{\tau}_i)=&-\max\{\gamma_{n,i}\phi_{n,i}(0)\check{W}_i^2(n),\hat{\rho}_i\tilde{\delta}_i(v_i)\}\\
  &+\max\{\gamma_{n,i}\phi_{n,i}(\hat{\tau}_i)\check{W}_i^2(n),\hat{\rho}_i\tilde{\delta}_i(v_i)\},
\end{split}\]
where $\tilde{W}_i(n)$ denotes $\tilde{W}_i(\bar{k}_i,n,\Theta_i^{n+1}, E_i^{1,n+1})$. Then, the event-triggering condition in (\ref{ET general}) can be implemented in a  conservative way:
\begin{equation}\label{implement case ii}
\begin{split}
f_{\eta}^i\le&-a_i\eta_i+\min_{n\in\{1,\dots,\hat{l}_i\}}\{\hat{\delta}_i(\check{v}_{i,n})\}+(1-\varepsilon_i)\tilde{\rho}_i\tilde{\delta}_i(\tilde{v}_i),\\
g_t^i(q)\le&\eta_i-\bar{\rho}_i\tilde{\delta}_i(v_i)+\bar{\rho}_i\max\{\tilde{\delta}_i(v_i),\varpi_i(\hat{\tau}_i)\tilde{\delta}_i(\tilde{v}_i)\}\\
&+\min_{n\in\{1,\dots,\hat{l}_i\}}\{\mathcal{W}_{i,n}^t(\hat{\tau}_i)\},\\
g_s^i(q)\le&\eta_i-\bar{\rho}_i\tilde{\delta}_i(v_i)+\bar{\rho}_i\max\{\tilde{\delta}_i(v_i),\varpi_i(\hat{\tau}_i)\tilde{\delta}_i(\tilde{v}_i)\}\\
&+\min_{n\in\{1,\dots,\hat{l}_i\}}\{\mathcal{W}_{i,n}^s(\hat{\tau}_i)\}.
\end{split}
\end{equation}

In addition, in the case of lacking both Items~(i)~and~(ii), one can implement a more conservative event-triggering condition by directly combining the ideas in (\ref{implement case i}) and (\ref{implement case ii})--the details are omitted.

\section{Simulations}
{In this section, we will illustrate the main results by two numerical examples. The first one involves non-globally Lipschitz dynamics and in the second, TOD scheduling protocols are considered for the case of two nodes sharing one network.}
\subsection{Example 1}
Consider the following nonlinear plant borrowed from \cite{Wangperiodic}:
\begin{equation}
  \nonumber
  \begin{cases}
    \dot{x}_{p,1}=d_1x_{p,1}^2-x_{p,1}^3+x_{p,2}+\hat{u}_1+w,\\
    \dot{x}_{p,2}=d_2x_{p,2}^2-x_{p,2}^3+x_{p,1}+\hat{u}_2+w,\\
     y_1=x_{p,1},\\
     y_2=x_{p,2},\\
  \end{cases}
\end{equation}
where $x_{p,1},x_{p,2}\in\mathbb{R}$ and the parameters $d_1=d_2=0.8$. A local static feedback controller, generating $\hat{u}_i=u_i=-2\hat{y}_i$, is implemented. Thus, from Remark \ref{rmk of zero controller}, one has $x=x_p$ and there are two communication channels that sample and transmit the output $y_i$ for updating $\hat{y}_i$ in the destination nodes for $i=1,2$. Consequently, we have $v=[v_1,v_2]^{\rm{T}}=[y_1,y_2]^{\rm{T}}\in\mathbb{R}^2$.

Then, some relevant functions in (\ref{flow closed loop}) are given as follows:
\[\begin{split}
  f(q,w)=\left[\begin{smallmatrix}
    d_1x_1^2-x_1^3+x_2-2(x_1+e_1)+w\\
    d_2x_2^2-x_2^3+x_1-2(x_2+e_2)+w
  \end{smallmatrix}\right], {\text{ }} g(q,w)=-f(q,w).
\end{split}\]
Since there is only one sensor node in each communication channel, the corresponding scheduling protocol is given by $h_{v_i}=0$ for $i=1,2$, which results in Assumption \ref{asm of protocol} with
\[W_i(\bar{k}_i,e_i)=\left|e_i\right|, {\text{ }} \lambda_i=0, {\text{ }}\lambda_{W,i}=1, {\text{ }} M_{p,i}=1,\]
and Assumption \ref{asm of fvi} with
\[M_{e,i}=2, {\text{ }} H_i(x,e,w)=\left|d_ix_i^2-x_i^3+x_{3-i}-2x_i\right|+\left|w\right|.\]
for $i=1,2$.

For properties on the flow dynamics of $x$, we start with the assumption of $V_{\rm{R}}$ in (\ref{derivative of V reference}), where we take $V_{\rm{R}}(x)=a^2\sum_{i=1}^2(b\frac{x_i^2}{2}+c\frac{x_i^4}{4})$ and $\delta_{{\rm{R}},i}(v_i)=\frac{y_i^2}{2}=\frac{v_i^2}{2}$ for $i=1,2$, with some $a,b,c>0$. According to the calculations in \cite{Wangperiodic}, we select $(a,b,c)=(1.7,3.93,2.9)$ and
\[J_{i}(x,e,w)=-2x_i^2+\left|x_i^3+x_ix_{3-i}\right|+\left|2d_ix_i^2\right|+\left|e_i^2\right|+\left|w\right|^2, \]
for $i=1,2$, which lead (\ref{derivative of V reference}) to be
  \begin{equation}
    \label{derivative of V1 simulation}
     \begin{split}
     \left\langle \nabla V_{\rm{R}}(x), f(q,w)\right\rangle\le& -0.01\left\|x\right\|^2+\sum_{i=1}^N\Big((8.36^2-0.01)W_i^2\\
     &-5\delta_{{\rm{R}},i}(v_i)-H_{i}^2(x,e,w)-{J}_{i}(x,e,w)\Big)\\
     &+81.6\left|w\right|^2.
  \end{split}
  \end{equation}
Since the coefficient before $\delta_{{\rm{R}},i}(v_i)$ is $-5<-1$, to have positive terms in $f_{\eta}^i$, it is not necessary to introduce $\epsilon_i$ as in Remark~\ref{rmk of derivative of V}. Then, according to Remark \ref{rmk of check}, one can obtain Assumption \ref{asm of V orignal}  and the resultant Assumptions~\ref{asm of we}--\ref{asm of V} from Propositions \ref{pps of we}--\ref{pps of V}, by choosing $V=V_{\rm{R}}$, $\delta_i=\delta_{{\rm{R}},i}$,
\[\begin{split}
 \hat{\delta}_i(s)=&\bar{\delta}_i(s)=2{(1-\check{\varepsilon}_i)}s^2, {\text{ }} s\in\mathbb{R}, \\
 \gamma_i^2=&8.36^2+\frac{2(1-\check{\varepsilon}_i)}{\check{\epsilon}_i},
\end{split}
\]
where $\check{\varepsilon}_1, \check{\varepsilon}_2\in(0,1],i=1,2,$ are user-specified parameters.

To compute $T_M^i$ from (\ref{condition TMI transmission}--\ref{condition TMI update}), we select $\tilde{\lambda}_1=\tilde{\lambda}_2=0.5$; and for any given $D_i$, we fix the initial values of (\ref{def of phi}) satisfying $\phi_{l_i,i}(0)=10$ for all $l_i\in\{0,\dots,D_i+1\}$. Table~\ref{table of TMI} gives the calculation results under different $\check{\varepsilon}:=\check{\varepsilon}_1=\check{\varepsilon}_2\in(0,1]$ and $D:=D_1=D_2\in\mathbb{Z}_{\ge0}$, which illustrate the analysis in Remark~\ref{rmk of delay assumption}.

\begin{table}[!hbtp]
\caption{Relationship between $T_M^i$ and $\epsilon$ as well as $D$.}
\begin{center}
\begin{tabular}{c|c|c|c|c}
\hline
\diagbox{$\check{\varepsilon}$}{$T_M^i$}{$D$} & $0$ & $ 1$ & $2$ & $ 3$\\
\hline
 $1$& $0.0109$s & $  0.0055$s & $0.0027$s & $ 0.0014$s \\
 $0.5$ & $0.0108$s & $ 0.0054$s & $0.0027$s & $ 0.0013$s \\
 $0.1$& $0.0098$s & $ 0.0049$s & $0.0025$s & $ 0.0012$s \\
 $0.01$ & $0.0058$s & $ 0.0029$s & $0.0014$s & $ 0.0007$s \\
 \hline
\end{tabular}
\end{center}
\label{table of TMI}
\end{table}

Figs. \ref{stateF}--\ref{ltF} show the simulation results in the case of $D=1$, where the other parameters are given by $x_1(0)=10,x_2(0)=-10,\check{\varepsilon}=0.5,a_1=a_2=0.01, \varepsilon_1=\varepsilon_2=0.5, \pi_1=\pi_2=0.99$, $w(t)=2\sin(20\pi t)$, and $\bar{\rho}_i,\hat{\rho}_i$ for $i\in\{1,2\}$ are selected as their maximum values in Theorem~\ref{thm of stability}. For simplicity, we assume that the two communication channels have full capacities in Section IV-C, which results in the implementation in (\ref{conditions on dynamics of eta}). The lower bound of inter-sampling periods is selected as $0.002$s. From Fig.~\ref{stateF}, the designed periodic event-triggered NCS is robust to the external disturbance. Figs.~\ref{sequencesF} and \ref{ltF} provide the evolutions of different time sequences and $l_i(t),i=1,2$, in the first $0.1$s. It can be observed that the ETs can actively discard unnecessary sampled outputs in the two independent communication channels which are subject to time-varying inter-sampling periods and large transmission delays. In Fig.~\ref{ltF}, at all the sampling instants, the left limit of $l_i(t),i=1,2$, is smaller than $D_i=1$, which illustrates Proposition \ref{pps of upper bound of l}. Moveover, the average inter-event times of two channels $\mathcal{C}_1$ and $\mathcal{C}_2$ are, respectively, $0.0159$s and $0.0149$s which are about two times larger than $T_M^i,i=1,2$.

\begin{figure}[!hbtp]
\centering
\includegraphics[width=9cm]{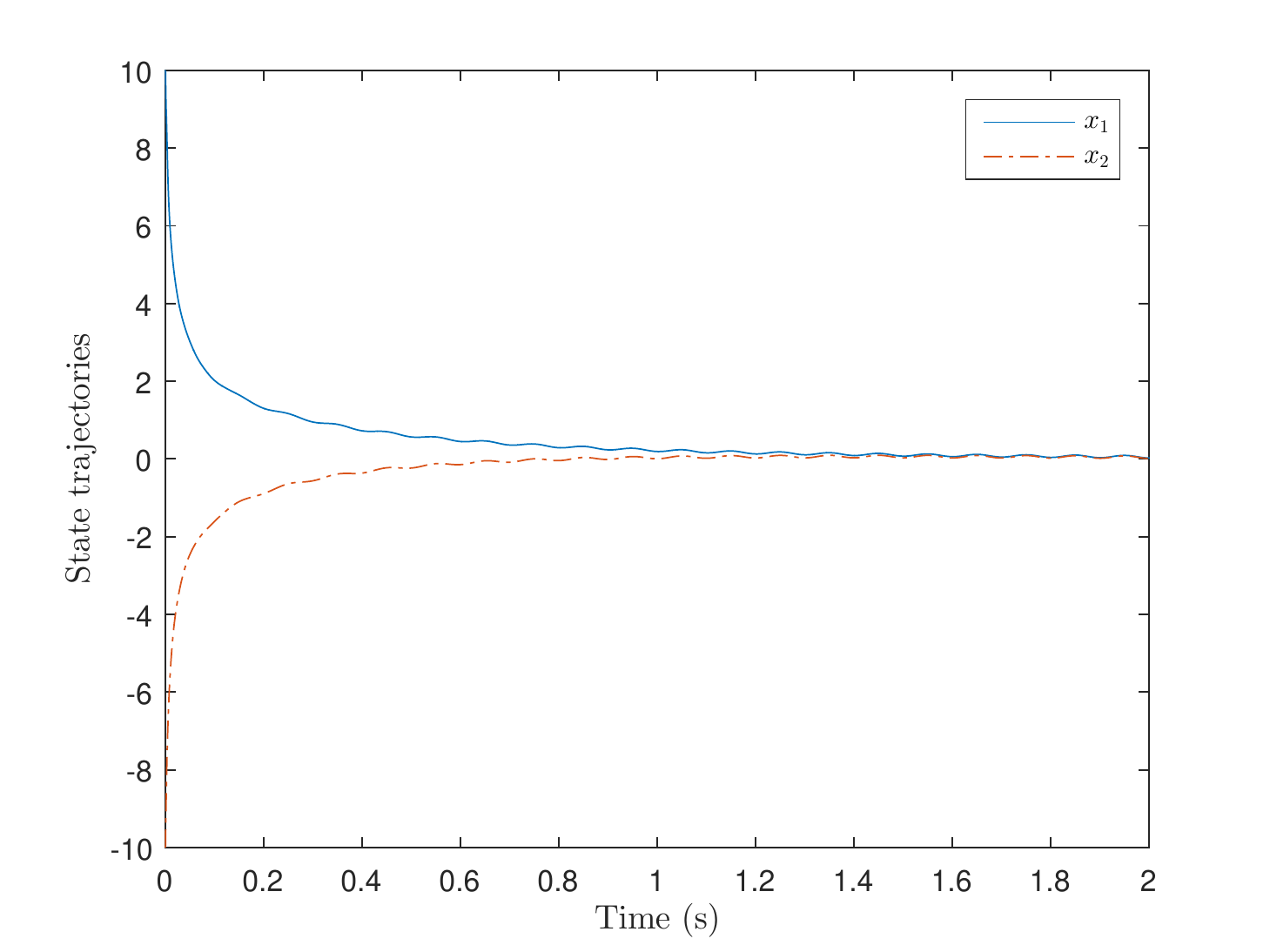}
\caption{State trajectories.}
\label{stateF}
\end{figure}

\begin{figure}[!hbtp]
\centering
\includegraphics[width=9cm]{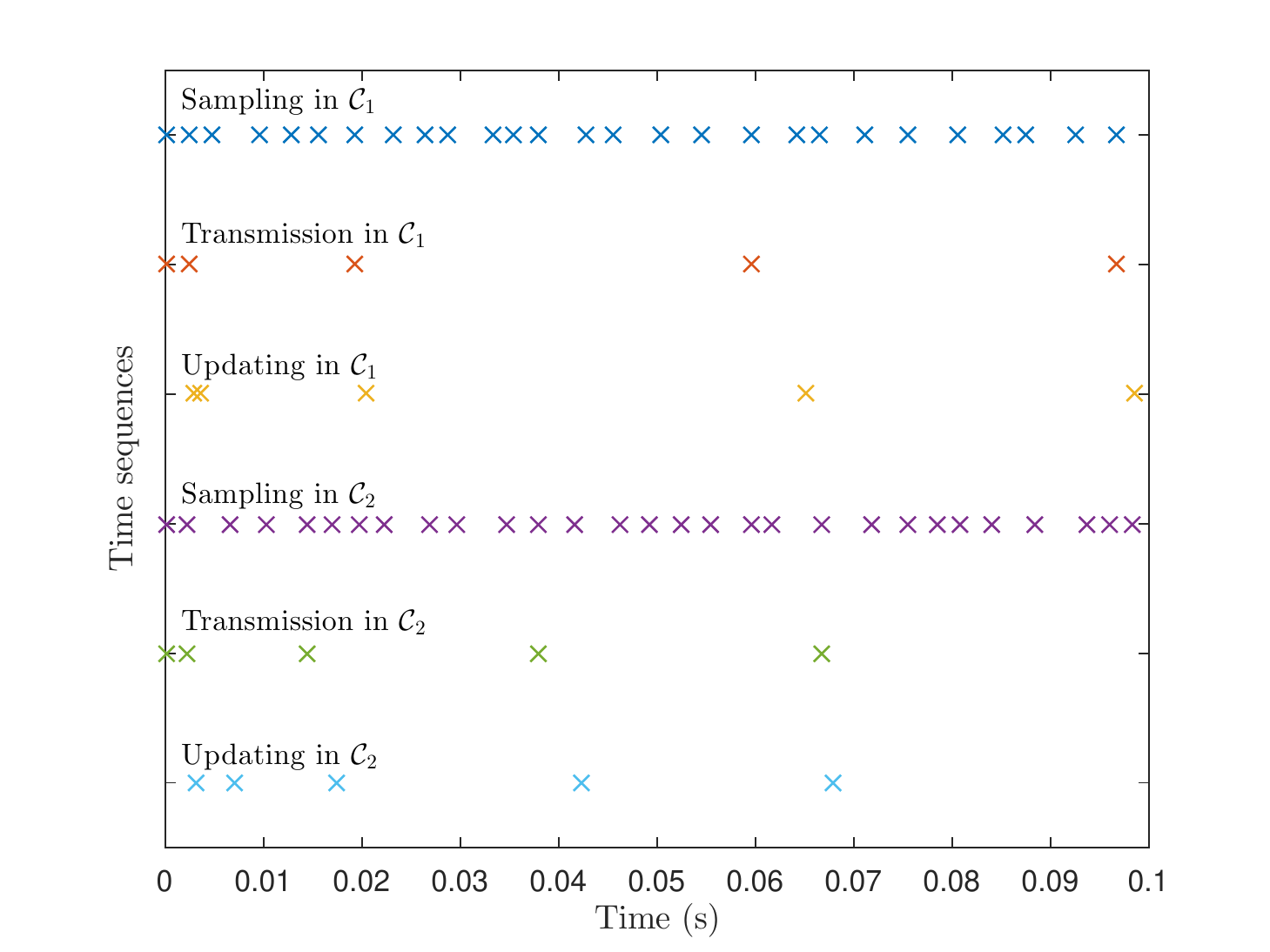}
\caption{Sampling, transmission, and updating time sequences in two communication channels.}
\label{sequencesF}
\end{figure}

\begin{figure}[!hbtp]
\centering
\includegraphics[width=9cm]{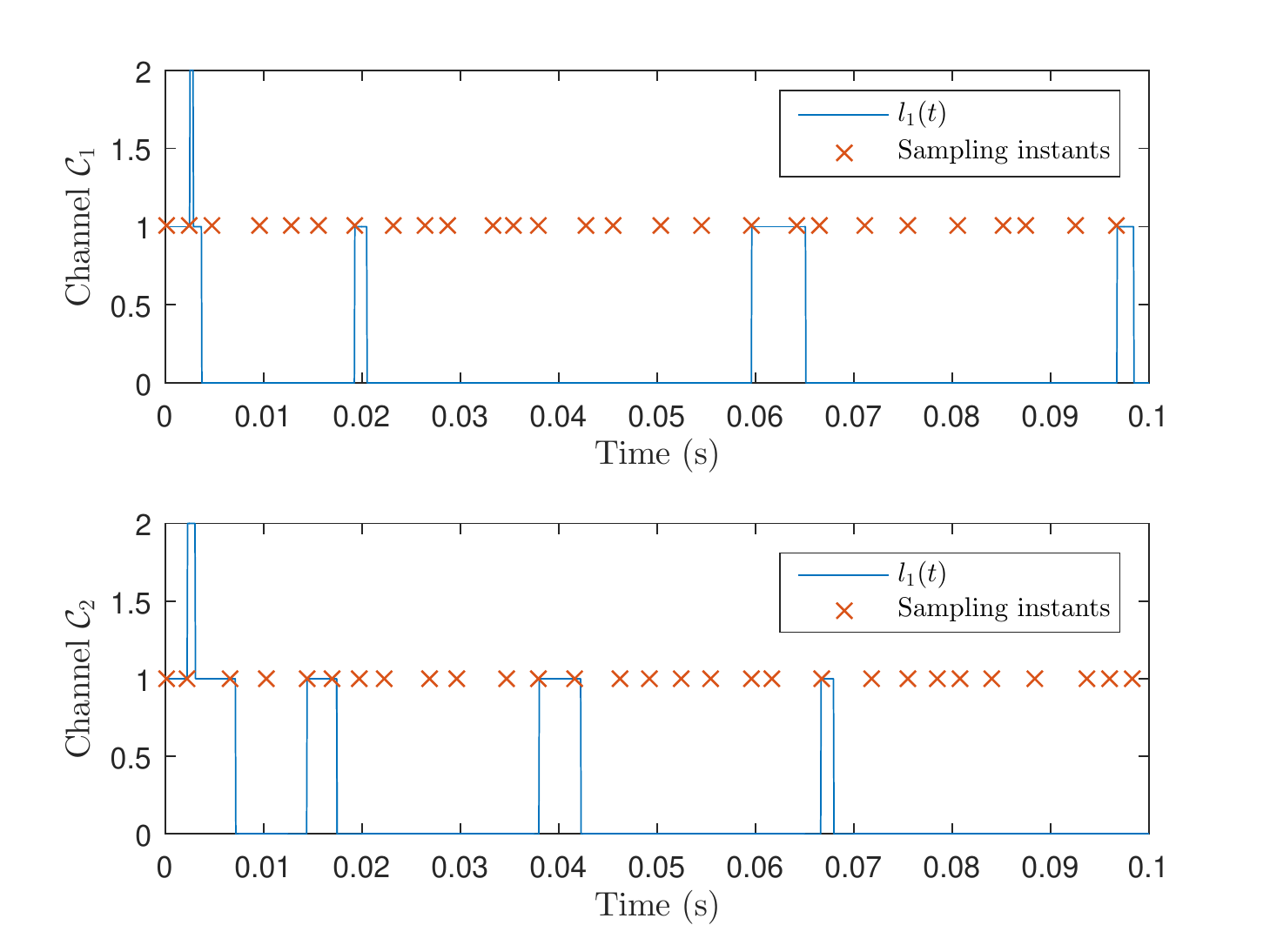}
\caption{The evolution of the variable $l_i(t),i=1,2$.}
\label{ltF}
\end{figure}

Finally, Table \ref{table of AIET} provides the simulations of the average inter-event times under different $\check{\varepsilon}$. The implementation is the same as in Figs. \ref{stateF}--\ref{ltF}. The average inter-event time in each case is obtained over $50$ simulations with random initial states. Table \ref{table of AIET} illustrates Remark \ref{rmk of check} by showing that allowing $\eta_i,i=1,2,$ to increase in flow can improve transmission performance although a large sampling frequency is required. This, in fact, also validates the advantages of the developed dynamic event-triggering conditions compared to the static ones in \cite{Wangperiodic}, where $\eta_i(t)=0$ for all $t\ge0$ and $i=1,2$, as discussed in Remark~\ref{rmk of static ET}.

\begin{table}[!hbtp]
\caption{Average inter-event times.}
\begin{center}
\begin{tabular}{c|c|c|c|c}
\hline
 $\check{\varepsilon}$ & $1$ & $ 0.5$ & $0.1$ & $ 0.01$\\
\hline
 $\mathcal{C}_1$& $0.0143$s & $ 0.0166$s & $0.0235$s & $  0.0313$s \\
 $\mathcal{C}_2$& $0.0144$s & $0.0166$s & $0.0234$s& $0.0312$s \\
 \hline
\end{tabular}
\end{center}
\label{table of AIET}
\end{table}

\subsection{Example 2}
{Consider the following nonlinear example of a single-link robot arm from \cite{Heijmans} with $x_{p}=[x_{p,1},x_{p,1}]^{\rm{T}}\in\mathbb{R}^2$:
\begin{equation}
  \nonumber
  \begin{cases}
    \dot{x}_{p,1}=x_{p,2},\\
    \dot{x}_{p,2}=-4.905\sin(x_{p,1})+2\hat{u},\\
     y=[x_{p,1},x_2]^{\rm{T}},
  \end{cases}
\end{equation}
and the following static controller:
\[\hat{u}=u=\frac{1}{2}(\sin(\hat{y}_{1})-\hat{y}_1-\hat{y}_2),\]
which leads to $v=v_1=x_p=x$. We assume that $y=[y_1,y_2]^{\rm{T}}$ is sensed by two nodes but transmitted trough only one network $\mathcal{C}_1$; and furthermore, a TOD scheduling protocol is considered {where $W_1(\bar{k}_1,e_1)=\left\|e_1\right\|$  and  $h_1(\bar{e}_1)=(I-\Phi(\bar{e}_1))\bar{e}_1$ with $\bar{e}_1=[\bar{e}_{1,1},\bar{e}_{1,2}]^{\rm{T}}$, $\Phi(\bar{e}_1):={\rm{diag}}\{\phi_1(\bar{e}_1),\phi_2(\bar{e}_1)\}$ and, for $j=1,2$,
\[\phi_j(\bar{e}_1):=\begin{cases}
  1, & {\text{if }} j=\min(\arg\max_{j}\left|\bar{e}_{1,j}\right|)\\
  0, & {\text{otherwise.}}
\end{cases}\]
More details on TOD scheduling protocol can be found in Example 2 of \cite{Nesic2004}.}}

{Then, based on the calculations in \cite{Heijmans}, the basic assumptions in Section~IV--B are given as follows:}
\begin{itemize}
  \item {For Assumption \ref{asm of protocol}: $W_1(\bar{k}_1,e_1)=\left\|e_1\right\|, \lambda_1=\sqrt{\frac{1}{2}}$ and  $\lambda_{W,1}=M_{p,1}=1$.}
  \item {For Assumption \ref{asm of fvi}: $H_1(x,e,0)=\left|x_2\right|+\left|x_1+x_2\right|$ and $M_{e,1}=8.351$.}
  \item {For Assumption \ref{asm of V orignal}: $V(x)= 12.2160x_1^2+4.2200x_1x_2+20.2120x_2^2$ and $\delta_1(v_1)=\left\|v_1\right\|^2$, leading to $J_1(x,e,0)=8.601\left\|e_1\right\|^2$ and
  \[\begin{split}
     \left\langle \nabla V(x), f(q,w)\right\rangle\le& -0.0116\left\|x\right\|^2+(\gamma_{1}^2-\bar{\varepsilon}_1)W_1^2\\
     &-\delta_1(v_1)-1.05H_{1}^2(x,e,0)\\
     &-0.0634\left\|\hat{v}_1\right\|^2-1.05{J}_1(x,e,0),
  \end{split}\]
  with $\gamma_1=46.1014$, $\bar{\varepsilon}_1=0.01$, and $\bar{\delta}_1(\hat{v}_1)=0.0634\left\|\hat{v}_1\right\|^2$.}
\end{itemize}

{For simulation, we consider the case of $D_1=1$ and $\tilde{\lambda}_1=0.8>\lambda_1$; then by solving Item 1) in Theorem \ref{thm of stability} with $\dot{\phi}_{l_1,1}=1.2$ for $l_1\in\{1,0\}$, one can choose $T_M^1=0.0016$s and $T_m^1=0.0005$s. Other parameters are given by $x_1(0)=-1.3,x_2(0)=2,\varepsilon_1=0.3, a_1=0.01,\pi_1=0.99$, and $\bar{\rho}_1,\hat{\rho}_1$  are selected as their maximum values. Fig. \ref{arm} provides the simulation results where the disturbance-free plant is asymptotically stabilized by our design dynamic PETC. The bottom sub-figure shows that only one node is scheduled in each transmission and most of unnecessary sampled signals will be discarded actively to save communication resources. Especially, the average inter-event time, $0.0045$, is about two times larger than $T_M^1$.}

\begin{figure}[!hbtp]
\centering
\includegraphics[width=9cm]{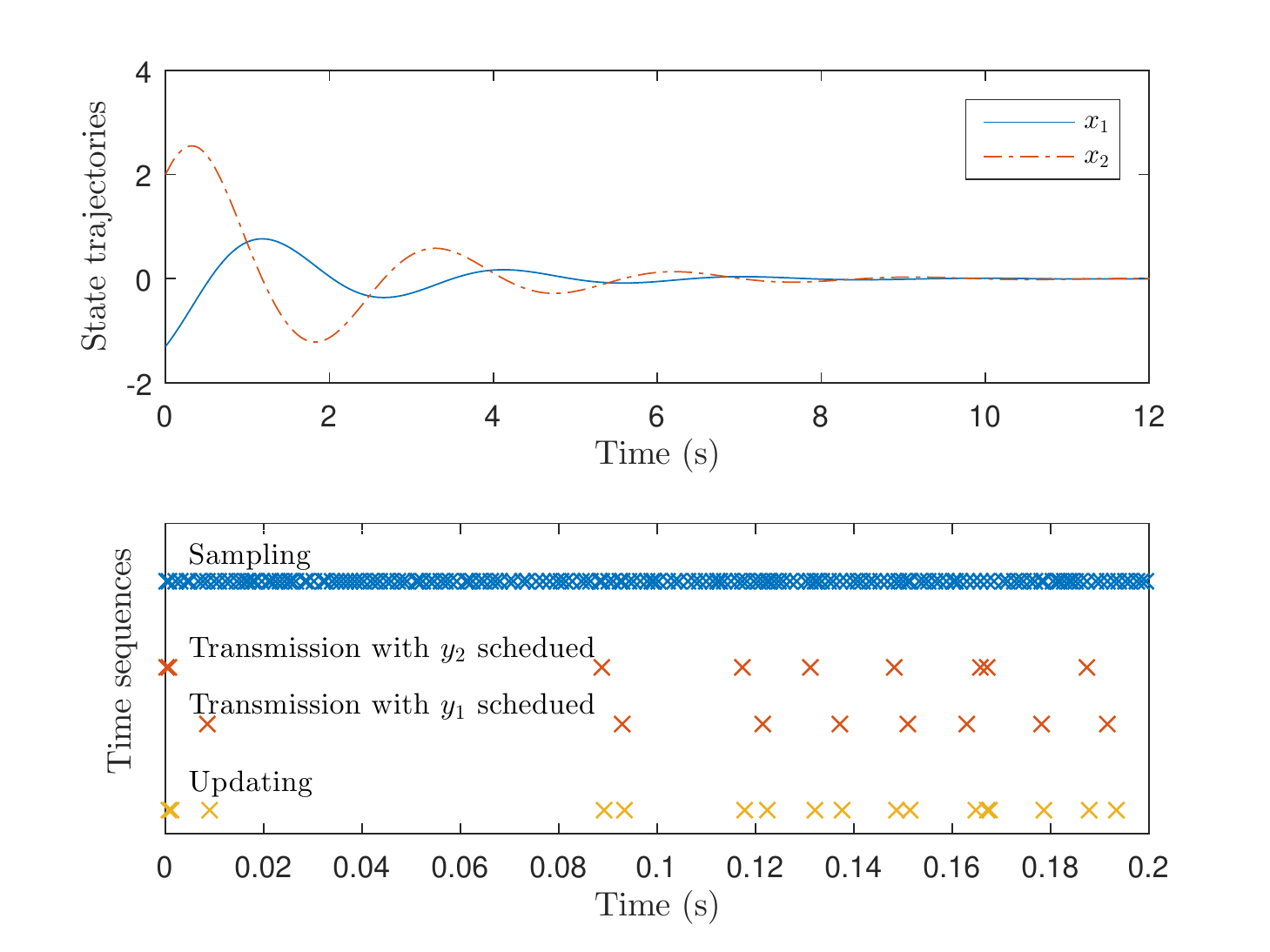}
\caption{{Simulation results in Example 2.}}
\label{arm}
\end{figure}

\section{Conclusions}
This paper studied periodic event-triggered networked control for nonlinear systems, where the plants and controllers were connected by multiple independent communication channels that were subject to lager transmission delays. Based on the assumption of MADNS, a new hybrid system approach was provided to model the closed-loop system with additional time-varying inter-sampling times, sensor node scheduling, and external disturbances. Then, by constructing new storage functions on the system state and updating errors, some inequalities were provided to characterize the relationship between MASP and MADNS. Moreover, to efficiently save limited communication resources, a new dynamic periodic event-triggered scheme was proposed, which included some existing static ones as special cases. From some assumptions provided by the emulation-based method, sufficient conditions on the design of dynamic ETC were given to ensure the input-to-state stability. Furthermore, according to different capacities of the equipment in communication channels, the implementation strategies of the designed dynamic event-triggered control were discussed. Finally, nonlinear examples was simulated to illustrate the feasibility and efficiency of the theoretical results. The simulation showed that the design of dynamic ETC might be improved by slightly increasing the sampling frequency. {Future work may include the investigation of disorder and distributed systems in the large-delay case.}

\section*{APPENDIX A: LEMMAS}
Some necessary technical lemmas are provided here.
\begin{lem}[\cite{Clarke}]\label{lem of clarke}
Consider two functions $U_1,U_2: \mathbb{R}^n \to \mathbb{R}$ that have well-defined
Clarke derivatives for $x \in \mathbb{R}^n$ and $v \in \mathbb{R}^n$. Introduce three sets
$\mathcal{A} := \{x : U_1(x)>U_2(x)\}$, $\mathcal{B} := \{x : U_1(x)<U_2(x)\}$, $\Gamma := \{x : U_1(x)=U_2(x)\}$.
Then for any $v \in \mathbb{R}^n$, the function $U(x) := \max\{U_1(x),U_2(x)\}$ satisfies
$U^{\circ}(x;v)=U^{\circ}_1(x;v)$ for $x \in \mathcal{A}$, $U^{\circ}(x;v)=U^{\circ}_2(x;v)$ for $x \in \mathcal{B}$,
and $U^{\circ}(x;v) \le \max\{U^{\circ}_1(x;v),U^{\circ}_2(x;v)\}$ for $x \in \Gamma$.
\end{lem}

\begin{lem}[\cite{Wangperiodic}]\label{lem of function K}
For any $N$ functions $\alpha_1,\dots,\alpha_N\in\mathcal{K}_{\infty}$, there exist $\mathcal{K}_{\infty}$--functions $\underline{\alpha}$ and $\bar{\alpha}$ such that $\underline{\alpha}(\sum_{i=1}^Ns_i)\le \sum_{i=1}^N \alpha_i(s_i)\le\bar{\alpha}(\sum_{i=1}^N s_i)$ holds for all $s_i\ge0$ and $i\in\bar{N}$.
\end{lem}

\section*{APPENDIX B: PROOFS}
\begin{IEEEproof}[Proof of Proposition \ref{pps of thetai}]
\textbf{Item 1):} At the sampling but not transmission instant $t\in\{s_j^i\}_{j\in\mathbb{Z}_{\ge0}}\backslash\{t_k^i\}_{k\in\mathbb{Z}_{\ge0}}$, neither $\bar{v}_i(t)$ nor $l_i(t)$ changes, thus there is no update for $\theta_i(t)$.

\textbf{Item 2):} At the updating instant $t=f_k^i$, the difference between transmission and updating numbers will decrease, i.e., $l_i(t^+)=l_i(t)-1$ while $\bar{k}_i(t)$ keeps constant. Thus, we have that for $\theta_i(t^+)$, only the first $l_i(t)-1$ blocks could be nonzero with $\theta_{i,j}(t^+)=\theta_{i,j+1}(t), j\in\{1,\dots, l_i(t)-1\}$. The other blocks are proved by considering that they are zero.

\textbf{Item 3):} At the transmission instant $t=t_k^i$, we have $l_i(t^+)=l_i(t)+1$ and $\bar{k}_i(t^+)=\bar{k}_i(t)+1=k+1$. Thus, in $\theta_i(t^+)$, the first $l_i(t)+1$ blocks could be nonzero and the first $l_i(t)$ blocks keep constant. For the $(l_i(t)+1)$-th block in $\theta_i(t^+)$, from (\ref{def of theta}), it follows that
\[\begin{split}
  \theta_{i,l_i(t)+1}(t^+)=&\bar{v}_i\left(t_{\bar{k}_i(t)}^{i+}\right)-\bar{v}_i\left(t_{\bar{k}_i(t)-1}^{i+}\right)\\
  =&\bar{v}_i\left(t_{\bar{k}_i(t)}^{i+}\right)-\bar{v}_i\left(t_{\bar{k}_i(t)-1}^{i+}\right)+\bar{v}_i\left(t_{\bar{k}_i(t)-2}^{i+}\right)\\
  &-\bar{v}_i\left(t_{\bar{k}_i(t)-2}^{i+}\right)\\
  =&\bar{v}_i\left(t_{\bar{k}_i(t)}^{i+}\right)-\bar{v}_i\left(t_{\bar{k}_i(t)-2}^{i+}\right)-\theta_{i,l_i(t)}(t)\\
  \vdots&\\
  =&\bar{v}_i\left(t_{\bar{k}_i(t)}^{i+}\right)-\bar{v}_i\left(t_{\bar{k}_i(t)-l_i(t)-1}^{i+}\right)-\sum_{j=1}^{l_i(t)}\theta_{i,j}(t)\\
  =&\bar{v}_i\left(t_{\bar{k}_i(t)}^{i+}\right)-\bar{v}_i\left(t_{\tilde{k}_i(t)-1}^{i+}\right)-\sum_{j=1}^{l_i(t)}\theta_{i,j}(t),
\end{split}\]
where the last equality is based on the definition of $l_i(t)$. Then from the updating rule in (\ref{update of hat vi}), it follows that
\[\begin{split}
  \theta_{i,l_i(t)+1}(t^+)=&\bar{v}_i\left(t_{\bar{k}_i(t)}^{i+}\right)-\hat{v}_i\left(f_{\tilde{k}_i(t)-1}^{i+}\right)-\sum_{j=1}^{l_i(t)}\theta_{i,j}(t)\\
  =&v_i(t)+h_{v_i}(\bar{k}(t),\bar{e}(t_{\bar{k}(t)}^{i}))-\hat{v}_i(t)-\sum_{j=1}^{l_i(t)}\theta_{i,j}(t)\\
  =&h_{v_i}(\bar{k}(t),\bar{e}_i(t_{\bar{k}(t)}^{i}))-e_i(t)-\sum_{j=1}^{l_i(t)}\theta_{i,j}(t),
\end{split}\]
where in the second equality, we utilize the facts of $t_{\bar{k}_i(t)}^{i}=t$ and $t\in(f_{\tilde{k}_i(t)-1}^{i},f_{\tilde{k}_i(t)}^{i}]$ based on the left continuity of $\bar{k}(t)$ and $\tilde{k}(t)$. In addition, following a similar process, we have
\[\begin{split}
  \bar{e}_i(t_{\bar{k}(t)}^{i})=&\bar{v}_i(t_{\bar{k}(t)}^{i})-v_i(t_{\bar{k}(t)}^{i})\\
  =&\bar{v}_i(t_{\bar{k}(t)-1}^{i+})-\bar{v}_i(t_{\bar{k}(t)}^{i+})+\bar{v}_i(t_{\bar{k}(t)}^{i+})-v_i(t_{\bar{k}(t)}^{i})\\
  \vdots&\\
  =&\sum_{j=1}^{l_i(t)}\theta_{i,j}(t)+\hat{v}_i\left(f_{\tilde{k}_i(t)-1}^{i+}\right)-v_i(t_{\bar{k}(t)}^{i})\\
  =&\sum_{j=1}^{l_i(t)}\theta_{i,j}(t)+e(t),
\end{split}
\]
which completes the proof.
\end{IEEEproof}

\begin{IEEEproof}[Proof of Theorem \ref{thm of stability}]
From the forms of $f_{\eta}^i$ and $g_t^i$ in the theorem and the dynamics of $\phi_{l_i,i}$ in (\ref{def of phi}) for all $i\in\bar{N}$ and $l_i\in\{0,\dots,D_i+1\}$, the condition in (\ref{condition TMI transmission}) ensures that $f_{\eta}^i(\cdot, 0)$ and $g_t^i(\cdot,\eta_i)$ are non-negative for all $\eta_i\ge0$. Thus, $\eta_i,i\in\bar{N}$, is non-negative scalars according to the analysis below (\ref{dynamics of etai}).

Consequently, for $q\in\mathcal{C}\cup\mathcal{D}$, we consider the following Lyapunov function:
\begin{equation}
  \label{def of U}
  \begin{split}
    U(q):=&\tilde{V}(x)+\sum_{i=1}^N \left(\bar{S}_i(q)+\hat{S}_i(q)+\eta_i\right)\\
    \bar{S}_i(q):=&\bar{\rho}_i\max\{\tilde{\delta}_i(v_i),\varpi_i(\hat{\tau}_i)\tilde{\delta}_i(\tilde{v}_i)\}\\
    \hat{S}_i(q):=&\max\{\gamma_{l_i,i}\phi_{l_i,i}(\hat{\tau}_i)\tilde{W}_i^2,\hat{\rho}_i\tilde{\delta}_i(v_i)\},
  \end{split}
\end{equation}
where, for $i\in\bar{N}$, the constants $\bar{\rho}_i,\hat{\rho}_i\ge0$ are defined in the theorem, and the variable $\varpi_i$ are defined in (\ref{def of varpi}). Meanwhile, introduce the following auxiliary  set:
\[\mathcal{S}_{\bar{\rho}}=\{q\in\mathcal{X}|x=0,e=0,\eta=0,\bar{\rho}_i\tilde{v}_i=0,i\in\bar{N}\},\]
and define $\bar{\rho}=[\bar{\rho}_1,\dots,\bar{\rho}_N]^{\rm{T}}$. From the definition of $\mathcal{S}$ in Definition \ref{definition of ISS}, it follows $\left\|q(t,\bar{j})\right\|_{\mathcal{S}}\le \left\|q(t,\bar{j})\right\|_{\mathcal{S}_{\bar{\rho}}}$ for $q\in\mathcal{X}$ and there exists a constant $\bar{c}>0$ such that $\left\|q(0,0)\right\|_{\mathcal{S}_{\bar{\rho}}}\le (1+\bar{c})\left\|q(0,0)\right\|_{\mathcal{S}}$ for all the initial state satisfying $\tilde{v}_i(0,0)=g_{v_i}(x(0,0))$, $i\in\bar{N}$. As a result, the input-to-state stability of $\mathcal{S}_{\bar{\rho}}$ is sufficient to the inequality in (\ref{inequality of ISS}). Hence, in the rest of this proof, we will show that $\mathcal{S}_{\bar{\rho}}$ is input-to-state stable for the closed-loop system in (\ref{closed loop system}--\ref{def of Gij3}).

Referring to Theorem 1 in \cite{Wangperiodic}, to prove the input-to-state stability of $\mathcal{S}_{\bar{\rho}}$, we only need to show the following properties: there exist $\mathcal{K}_{\infty}$--functions $\underline{\beta}_U,\bar{\beta}_U,\alpha_U,\alpha_F$ such that
\begin{enumerate}
  \item[I.] $U$ is locally Lipschitz in $(x,e,\tilde{v},\hat{\tau})$, and for all $q\in\mathcal{C}\cup\mathcal{D}$ satisfying $\theta_{i,j}=0$ with $i\in\bar{N}$ and $j=l_i+1,\dots,D_i+1$,
  \[\underline{\beta}_U(\left\|q\right\|_{\mathcal{S}_{\bar{\rho}}})\le U(q) \le\bar{\beta}_U(\left\|q\right\|_{\mathcal{S}_{\bar{\rho}}});\]
  \item[II.] for all $q\in\mathcal{C}$ and $w\in\mathbb{R}^{n_w}$,
  \[U^{\circ}(q;F(q,w))\le -\alpha_U(U(q))+\alpha_F(\left\|w\right\|);\]
  \item[III.] for all $q\in\mathcal{D}, w\in\mathbb{R}^{n_w}$, and $g(q)\in G(q)$,
  \[U(g(q))\le U(q).\]
\end{enumerate}

\textbf{First, we consider the proof of Item I}. From Assumptions~\ref{asm of we}--\ref{asm of V} and the definitions of $\phi_{l_i,i}$ and $\varpi_i$ for $i\in\bar{N}$ and $l_i\in\{0,\dots, D_i+1\}$, $U$ satisfies the local Lipschitz property. Then according to Theorem 1 in \cite{Wangperiodic}, the continuity and positive definiteness of $\tilde{\delta}_i$ and the fact that the continuously differentiable function $g_{v_i}$ satisfies $v_i=g_{v_i}(x)$ and $g_{v_i}(0)=0$ imply that there exists a $\mathcal{K}$--function $\alpha_{\tilde{\delta}_i}$ satisfying $\tilde{\delta}_i(v_i)\le\alpha_{\tilde{\delta}_i}(\left\|[x^{\rm{T}},e^{\rm{T}}]^{\rm{T}}\right\|)$. Meanwhile, from the boundary conditions on $\phi_{l_i,i}$ in the theorem and on $\varpi_i$ in (\ref{def of varpi}) for $i\in\bar{N}$ and $l_i\in\{0,\dots,D_i+1\}$, we have that the positive constants $\underline{\phi}_i,\bar{\phi}_i$ and $\underline{\varpi}_i>0$ in (\ref{bound of phi}) are well defined.

Recall the definitions of constants, $\underline{\gamma}_i$ and $\bar{\gamma}_i$ in (\ref{bound of gamma}). Then, from (\ref{positive of wi}), (\ref{positive of V}), and (\ref{def of U}), it follows that
\[
\begin{split}
  U(q)\ge& \underline{\beta}_{\tilde{V}}(\left\|x\right\|)+\sum_{i=1}^N\left(\eta_i+\bar{\rho}_i\underline{\varpi}_i\underline{\beta}_{\tilde{\delta}_i}(\left\|\tilde{v}_i\right\|)\right)\\
  &+\sum_{i=1}^N\underline{\gamma}_i\underline{\phi}_i \underline{\beta}_{\tilde{W}_i}^2(\left\|[e_i^{\rm{T}},\theta_i^{\rm{T}}]^{\rm{T}}\right\|),\\
  U(q)\le& \bar{\beta}_{\tilde{V}}(\left\|x\right\|)+\sum_{i=1}^N\left(\bar{\rho}_i\bar{\beta}_{\tilde{\delta}_i}(\left\|\tilde{v}_i\right\|)+
  \bar{\rho}_i\alpha_{\tilde{\delta}_i}(\left\|[x^{\rm{T}},e^{\rm{T}}]^{\rm{T}}\right\|)\right)\\
  &+\sum_{i=1}^N(\eta_i+\bar{\gamma}_i\bar{\phi}_i \bar{\beta}_{\tilde{W}_i}^2(\left\|[e_i^{\rm{T}},\theta_i^{\rm{T}}]^{\rm{T}}\right\|))\\
  &+\sum_{i=1}^N\hat{\rho}_i\alpha_{\tilde{\delta}_i}(\left\|[x^{\rm{T}},e^{\rm{T}}]^{\rm{T}}\right\|),
\end{split}\]
which proves Item I by using Lemma \ref{lem of function K}.

\textbf{Next, we consider the proof of Item II}. From the property of Clarke derivative in Lemma \ref{lem of clarke}, for each $i\in\bar{N}$, we distinguish the following nine cases: Case $(1,m)$  $\times$ Case $(2,n)$ with $m,n={1,2,3}$, where Case $(1,1)$: $\tilde{\delta}_i(v_i)>\varpi_i(\hat{\tau}_i)\tilde{\delta}_i(\tilde{v}_i)$;  Case $(1,2)$: $\tilde{\delta}_i(v_i)<\varpi_i(\hat{\tau}_i)\tilde{\delta}_i(\tilde{v}_i)$; Case $(1,3)$: $\tilde{\delta}_i(v_i)=\varpi_i(\hat{\tau}_i)\tilde{\delta}_i(\tilde{v}_i)$; Case $(2,1)$: $\gamma_{l_i,i}\phi_{l_i,i}(\hat{\tau}_i)\tilde{W}_i^2>\hat{\rho}_i\tilde{\delta}_i(v_i)$; Case $(2,2)$: $\gamma_{l_i,i}\phi_{l_i,i}(\hat{\tau}_i)\tilde{W}_i^2<\hat{\rho}_i\tilde{\delta}_i(v_i)$; and Case $(2,3)$: $\gamma_{l_i,i}\phi_{l_i,i}(\hat{\tau}_i)\tilde{W}_i^2=\hat{\rho}_i\tilde{\delta}_i(v_i)$. Denote by $\mathcal{N}_{3(m-1)+n}$ all the communication channels that belong to the Case $(1,m)$ and Case $(2,n)$ simultaneously. Thus, we have $\mathcal{N}_m\in\bar{N},\cup_{m=1}^9\mathcal{N}_j=\bar{N}$ and $\mathcal{N}_{m_1}\cap\mathcal{N}_{m_2} =\emptyset$ for $m,m_1,m_2\in\{1,\dots,9\}$ and $m_1\neq m_2$. {Note that when $\hat{\rho}_i=0$, Case (2,2) is impossible; then the corresponding analysis is vacuously true.}

From (\ref{dynamics of etai}) and (\ref{flow of V}), we have
\begin{equation}
  \label{derivative of U original}
  \begin{split}
     U^{\circ}(q;F(q,w))\le& -\alpha_{\tilde{V}}(\left\|x\right\|)+\alpha_w(\left\|w\right\|)-\sum_{i=1}^N\sigma_{l_i,i}(\tilde{W}_i)\\
    &+\sum_{i=1}^N\left(\bar{S}_i^{\circ}(q;F(q,w))+\hat{S}_i^{\circ}(q;F(q,w))\right)\\
    &+\sum_{i=1}^NZ_i(q,w),\\
    Z_i(q,w):=&\gamma_{l_i,i}^2\tilde{W}_i^2-\tilde{\delta}_i(v_i)-(1+\epsilon_{l_i,i})(H_{l_i,i}^2+\tilde{J}_i)\\
    &-\hat{\delta}_i(\hat{v}_i)+f_{\eta}^i(q),
  \end{split}
\end{equation}
where for simplicity, we omit the arguments of $H_{l_i,i}$ and $\tilde{J}_i$. In the following, we consider $\bar{S}_i^{\circ}(q;F(q,w))+\hat{S}_i^{\circ}(q;F(q,w))$ for $i$ in different cases, i.e., $i\in\mathcal{N}_{m},m\in\{1,\dots,9\}$. Define $S_i(q):=\bar{S}_i(q)+\hat{S}_i(q)$.

For $i\in\mathcal{N}_1$, we have $\tilde{\delta}_i(v_i)>\varpi_i(\hat{\tau}_i)\tilde{\delta}_i(\tilde{v}_i)$ and
\[S_i(q)=\bar{\rho}_i\tilde{\delta}_i(v_i)+\gamma_{l_i,i}\phi_{l_i,i}(\hat{\tau}_i)\tilde{W}_i^2,\]
which, from (\ref{flow of wi}), (\ref{flow of deltai}) and (\ref{def of phi}), leads to
\[\begin{split}
  {S}_i^{\circ}(q;F(q,w))\le& \bar{\rho}_i(H_{l_i,i}^2+\tilde{J}_i)+\gamma_{l_i,i}\phi_{l_i,i}\tilde{W}_i(L_{l_i,i}\tilde{W}_i+H_{l_i,i})\\
  &-\gamma_{l_i,i}(2L_{l_i,i}\phi_{l_i,i}+\gamma_{l_i,i}(\phi_{l_i,i}^2+1))\tilde{W}_i^2\\
  =&\bar{\rho}_i(H_{l_i,i}^2+\tilde{J}_i)+\gamma_{l_i,i}\phi_{l_i,i}W_iH_{l_i,i}\\
  &-\gamma_{l_i,i}^2(\phi_{l_i,i}^2+1)\tilde{W}_i^2,
\end{split}\]
where the argument of $\phi_{l_i,i}$ is omitted. Since $\bar{\rho}_i\le\min_{l_i\in\{0,\dots,D_i+1\}}\{\epsilon_{l_i,i}\}$, one has
\begin{equation}
\nonumber
\begin{split}
    Z_i(q,w)+{S}_i^{\circ}(q;F(q,w))\le -\tilde{\delta}_i(v_i)-\hat{\delta}_i(\hat{v}_i)+f_{\eta}^i(q).
\end{split}
\end{equation}
The dynamics in (\ref{def of varpi}) implies $\varpi_i(\hat{\tau}_i)\ge \pi_i$, which yields
\begin{equation}\label{case 1 inequality}
\begin{split}
    Z_i(q,w)+{S}_i^{\circ}(q;F(q,w))\le -\pi_i\tilde{\delta}_i(\tilde{v}_i)-\hat{\delta}_i(\hat{v}_i)+f_{\eta}^i(q).
\end{split}
\end{equation}

For $i\in\mathcal{N}_2$, we have $\gamma_{l_i,i}\phi_{l_i,i}(\hat{\tau}_i)\tilde{W}_i^2<\hat{\rho}_i\tilde{\delta}_i(v_i)$ and
\[S_i(q)=\bar{\rho}_i\tilde{\delta}_i(v_i)+\hat{\rho}_i\tilde{\delta}_i(v_i).\]
From $\hat{\rho}_i\le\frac{1}{2}\min\{1,\underline{\phi}_i/\bar{\gamma}_i\}$, it follows
\[\gamma_{l_i,i}^2\tilde{W}_i^2\le \frac{\bar{\gamma}_i\hat{\rho}_i\tilde{\delta}_i(v_i)}{\phi_{l_i,i}(\hat{\tau}_i)}\le\frac{\bar{\gamma}_i\hat{\rho}_i\tilde{\delta}_i(v_i)}{\underline{\phi}_i}\le \frac{1}{2}\tilde{\delta}_i(v_i), \]
for all $l_i\in\{0,\dots,D_i+1\}$. Consequently, we have
\begin{equation}\label{case 2 inequality}
\begin{split}
    Z_i(q,w)+{S}_i^{\circ}(q;F(q,w))\le -\frac{\pi_i}{2}\tilde{\delta}_i(\tilde{v}_i)-\hat{\delta}_i(\hat{v}_i)+f_{\eta}^i(q).
\end{split}
\end{equation}

For $i\in\mathcal{N}_4$, we have $\tilde{\delta}_i(v_i)<\varpi_i(\hat{\tau}_i)\tilde{\delta}_i(\tilde{v}_i)$ and
\[S_i(q)=\bar{\rho}_i\varpi_i(\hat{\tau}_i)\tilde{\delta}_i(\tilde{v}_i)+\gamma_{l_i,i}\phi_{l_i,i}(\hat{\tau}_i)\tilde{W}_i^2,\]
which implies from (\ref{def of varpi}) that
\[\begin{split}
  {S}_i^{\circ}(q;F(q,w))\le& -\frac{\bar{\rho}_i(1-\pi_i)}{T_M^i}\tilde{\delta}_i(\tilde{v}_i)+\gamma_{l_i,i}\phi_{l_i,i}W_iH_{l_i,i}\\
  &-\gamma_{l_i,i}^2(\phi_{l_i,i}^2+1)\tilde{W}_i^2,
\end{split}\]
and subsequently,
\begin{equation}\label{case 4 inequality}
\begin{split}
    Z_i(q,w)+{S}_i^{\circ}(q;F(q,w))\le& -\frac{\bar{\rho}_i(1-\pi_i)}{T_M^i}\tilde{\delta}_i(\tilde{v}_i)-\tilde{\delta}_i(v_i)\\
    &-\hat{\delta}_i(\hat{v}_i)+f_{\eta}^i(q).
\end{split}
\end{equation}

For $i\in\mathcal{N}_5$, we have
\[S_i(q)=\bar{\rho}_i\varpi_i(\hat{\tau}_i)\tilde{\delta}_i(\tilde{v}_i)+\hat{\rho}_i\tilde{\delta}_i(v_i).\]
Following a similar process of (\ref{case 2 inequality}) and (\ref{case 4 inequality}), we have
\begin{equation}\label{case 5 inequality}
\begin{split}
    Z_i(q,w)+{S}_i^{\circ}(q;F(q,w))\le& -\frac{\bar{\rho}_i(1-\pi_i)}{T_M^i}\tilde{\delta}_i(\tilde{v}_i)-\frac{1}{2}\tilde{\delta}_i(v_i)\\
    &-\hat{\delta}_i(\hat{v}_i)+f_{\eta}^i(q).
\end{split}
\end{equation}

From (\ref{conditions on dynamics of eta}a), the flow of $\eta_i$ satisfies
\[\begin{split}
  f_{\eta}^i(q)\le& -a_i\eta_i+\hat{\delta}_i(\hat{v}_i)+(1-\varepsilon_i)\tilde{\rho}_i\tilde{\delta}_i(\tilde{v}_i),\\
  \tilde{\rho}_i:=&\min\left\{\frac{\bar{\rho}_i(1-\pi_i)}{T_M^i},\frac{\pi_i}{2}\right\},
\end{split}\]
with some $\varepsilon_i\in(0,1)$ and $a_i>0$. Thus, combining (\ref{case 1 inequality}--\ref{case 5 inequality}) yields
\begin{equation}\label{union case}
\begin{split}
  &\sum_{i\in\mathcal{N}_{\rm{SI}}}(Z_i(q,w)+{S}_i^{\circ}(q;F(q,w)))\le \sum_{i\in\mathcal{N}_{\rm{SI}}} -(\varepsilon_i\tilde{\rho}_i\tilde{\delta}_i(\tilde{v}_i)+a_i\eta_i),
\end{split}
\end{equation}
where $\mathcal{N}_{\rm{SI}}:=\mathcal{N}_1\cup\mathcal{N}_2\cup\mathcal{N}_4\cup\mathcal{N}_5$ and the subscript $\rm{SI}$ stands for ``strict inequalities''.

For the rest cases: $i\in\bar{N}\backslash\mathcal{N}_{\rm{SI}}$, from Lemma \ref{lem of clarke} and (\ref{union case}), one also has
\[Z_i(q,w)+{S}_i^{\circ}(q;F(q,w))\le -(\varepsilon_i\tilde{\rho}_i\tilde{\delta}_i(\tilde{v}_i)+a_i\eta_i).\]
Therefore, applying all the cases to (\ref{derivative of U original}) leads to
\begin{equation}
  \label{flow of U}
  \begin{split}
    U^{\circ}(q;F(q,w))\le& -\sum_{i=1}^N \left(\sigma_{l_i,i}(\tilde{W}_i)+\varepsilon_i\tilde{\rho}_i\tilde{\delta}_i(\tilde{v}_i)+a_i\eta_i\right)\\
    &-\alpha_{\tilde{V}}(\left\|x\right\|)+\alpha_w(\left\|w\right\|).
  \end{split}
\end{equation}
Thus, from Assumptions \ref{asm of we} and \ref{asm of V}, there exists a $\mathcal{K}_{\infty}$ function $\alpha_U$ satisfying Item II with $\alpha_F=\alpha_w$.

\textbf{Finally, we consider the proof of Item III}. Since $\tilde{V}(x)$ is continuous, only the terms $\left(S_i(q)+\eta_i\right)$ needs consideration, where $S_i(q)$ is defined below (\ref{derivative of U original}). For each communication channel $\mathcal{C}_i,i\in\bar{N}$, we distinguish four cases based on (\ref{jump closed loop}). Case 1: $\hat{m}_i=1$ and $g_s^i<0$; Case 2: $\hat{m}_i=1$ and $g_s^i>0$; Case 3: $\hat{m}_i=-1$; and Case 4: $\hat{m}_i=1$ and $g_s^i=0$.

In Case 1, the current instant is a transmission instant with $l_i\le D_i$ and $G_i(q)=G_i^1(q)$, where $\hat{\tau}_i^+=0,\tilde{v}_i^+=v_i, \bar{k}_i^+=\bar{k}_i+1, l_i^+=l_i+1, e_i^+=e_i$ and $\eta_i^+=g_t^i(q)$. Thus, for $g_i\in G_i^1(q)$, from (\ref{jump of wi}a), we have
\[\begin{split}
  S_i(g_i)-S_i(q)\le&\bar{\rho}_i\tilde{\delta}_i(v_i)-\bar{\rho}_i\max\{\tilde{\delta}_i(v_i),\varpi_i(\hat{\tau}_i)\tilde{\delta}_i(\tilde{v}_i)\}\\
  &+\{\gamma_{l_i+1,i}\phi_{l_i+1,i}(0)\tilde{\lambda}_i^2\tilde{W}_i^2,\hat{\rho}_i\tilde{\delta}_i(v_i)\}\\
  &-\max\{\gamma_{l_i,i}\phi_{l_i,i}(\hat{\tau}_i)\tilde{W}_i^2,\hat{\rho}_i\tilde{\delta}_i(v_i)\}
\end{split}\]
due to $\varpi_{i}(0)=1$. Thus, according to the inequality of $g_t^i$ in (\ref{conditions on dynamics of eta}b), one has
\[S_i(g_i)+g_t^i(q)\le S_i(q)+\eta_i.\]

In Case 2, we have $\hat{\tau}_i^+=0, \tilde{v}_i^+=v_i, l_i^+=l_i,\bar{k}_i^+=\bar{k}_i, e_i^+=e_i$ and $\eta_i^+=g_s^i(q)>0$, which implies
\[\begin{split}
  S_i(g_i)-S_i(q)\le&\bar{\rho}_i\tilde{\delta}_i(v_i)-\bar{\rho}_i\max\{\tilde{\delta}_i(v_i),\varpi_i(\hat{\tau}_i)\tilde{\delta}_i(\tilde{v}_i)\}\\
  &+\max\{\gamma_{l_i,i}\phi_{l_i,i}(0)\tilde{W}_i^2,\hat{\rho}_i\tilde{\delta}_i(v_i)\}\\
  &-\max\{\gamma_{l_i,i}\phi_{l_i,i}(\hat{\tau}_i)\tilde{W}_i^2,\hat{\rho}_i\tilde{\delta}_i(v_i)\},
\end{split}\]
with $g_i= G_i^2(q)$. Thus, the condition of $g_s^i$ in (\ref{conditions on dynamics of eta}c) ensures
\[S_i(g_i)+g_t^i(q)\le S_i(q)+\eta_i.\]

In Case 3, the jump is due to the update of $\hat{v}_i$ at the destination node, which yields that $\hat{\tau}_i^+=\tau_i, \tilde{v}_i^+=\tilde{v}_i, l_i^+=l_i-1\ge 0,\bar{k}_i^+=\bar{k}_i, e_i^+=e_i+\theta_{i,1}$ and $\eta_i^+=\eta_i$. Thus, (\ref{jump of wi}b) and the condition in (\ref{condition TMI update}) ensure that for $g_i\in G_i^3(q)$,
\[S_i(g_i)+\eta_i\le S_i(q)+\eta_i.\]

In Case 4, from the analysis in Cases 1 and 2, one also has for $g_i\in \{G_i^1(q),G_i^2(q)\}$,
\[S_i(g_i)+g_t^i(q)\le S_i(q)+\eta_i.\]
Consequently, Item III can be proved by combining the analysis on Cases 1-4.

Therefore, following the same line of the proof for Theorem~1 in \cite{Wangperiodic}, one can conclude the input-to-state stability of the set $\mathcal{S}$ with some $\mathcal{KL}$--function $\beta$ and $\mathcal{K}_{\infty}$--function directly from Items I-III.
\end{IEEEproof}

\begin{IEEEproof}[Proof of Proposition \ref{pps of we}]
In Assumption \ref{asm of protocol}, the continuity of $W_i(\bar{k}_i,\cdot)$ for all fixed $\bar{k}_i\in\mathbb{Z}_{\ge0}$ and Item 4) ensure that $W_i(\bar{k}_i,\cdot)$ is globally Lipschitz for all fixed $\bar{k}_i\in\mathbb{Z}_{\ge0}$. Thus, the ``max'' form in (\ref{def of tilde wi}) ensures that $\tilde{W}_i(\bar{k}_i,l_i,\cdot,\cdot)$ is (at least) locally Lipschitz for all $\bar{k}_i\in\mathbb{Z}_{\ge0}$ and $l_i\in\{0,\dots,D_i+1\}$.

Then, we prove (\ref{positive of wi}). According to the form of $\theta_i$ defined in (\ref{def of theta}), for any given $l_i\in\{0,\dots,D_i+1\}$, we have
\[\left\|\theta_i\right\|=\left\|\left[\theta_{i,1}^{\rm{T}},\dots,\theta_{i,l_i}^{\rm{T}}\right]^{\rm{T}}\right\|,\]
which implies that there exist positive constants $\bar{m}_{i,l_i}\ge\underline{m}_{i,l_i}>0$ such that
\[\underline{m}_{i,l_i}\left\|[e_i^{\rm{T}},\theta_i^{\rm{T}}]^{\rm{T}}\right\|\le\sum_{n=0}^{l_i}\left\|e_i+\sum_{j=1}^{n}\theta_{i,j}\right\| \le \bar{m}_{i,l_i}\left\|[e_i^{\rm{T}},\theta_i^{\rm{T}}]^{\rm{T}}\right\|,\]
for all $[\bar{k}_i,l_i,e_i^{\rm{T}},\theta_i^{\rm{T}}]^{\rm{T}}\in \tilde{\mathcal{X}}_i$.

Hence, from Item 1) in Assumption \ref{asm of protocol} and Lemma \ref{lem of function K}, there exist $\mathcal{K}_{\infty}$--functions $\underline{\alpha}_{i,l_i}^{\tilde{W}}$ and $\bar{\alpha}_{i,l_i}^{\tilde{W}}$ satisfying
\[\begin{split}
\tilde{W}_i(\bar{k}_i,l_i,\theta_i, e_i)\ge&   \frac{\tilde{\lambda}_i^{l_i}}{\lambda_{W,i}^{l_i}}\underline{\alpha}_{i,l_i}^{\tilde{W}}(\underline{m}_{i,l_i}\left\|[e_i^{\rm{T}},\theta_i^{\rm{T}}]^{\rm{T}}\right\|),   \\
\tilde{W}_i(\bar{k}_i,l_i,\theta_i, e_i)\le&     \bar{\alpha}_{i,l_i}^{\tilde{W}}(\bar{m}_{i,l_i}\left\|[e_i^{\rm{T}},\theta_i^{\rm{T}}]^{\rm{T}}\right\|).
\end{split}\]
Consequently, (\ref{positive of wi}) is proved by selecting the $\mathcal{K}_{\infty}$--function $\underline{\beta}_{W_i}(\cdot)$ ($\bar{\beta}_{W_i}(\cdot)$) as the minimum (maximum) of $\frac{\tilde{\lambda}_i^{l_i}}{\lambda_{W,i}^{l_i}}\underline{\alpha}_{i,l_i}^{\tilde{W}}(\underline{m}_{i,l_i}\cdot)$ ($\bar{\alpha}_{i,l_i}^{\tilde{W}}(\bar{m}_{i,l_i}\cdot)$) over all $l_i\in\{0,\dots,D_i+1\}$.

Next, we consider the jump behavior of $\tilde{W}_i$ in (\ref{jump of wi}). If a transmission occurs at the current instant, then with the form of $\tilde{W}_i$ in (\ref{def of tilde wi}), we have
\begin{equation}
\nonumber
  \begin{split}
  \tilde{W}_i&(\bar{k}_i+1,l_i+1,[\theta_{i,1}^{\rm{T}},\dots,\theta_{i,l_i}^{\rm{T}},h_{v_i}^{\rm{T}}-e_i^{\rm{T}}-\sum_{j=1_i}^{l}\theta_{i,j}^{\rm{T}},0^{\rm{T}}]^{\rm{T}}, e_i)\\
  &= \max\Big\{\frac{\tilde{\lambda}_i^{l_i+1}}{\lambda_{W,i}^{l_i+1}}W_i(\bar{k}_i+1,e_i),\frac{\tilde{\lambda}_i^{l_i}}{\lambda_{W,i}^{l_i}}W_i(\bar{k}_i+1,e_i+\theta_{i,1}),\\
   &{\text{\quad\quad}}\dots,\frac{\tilde{\lambda}_i}{\lambda_{W,i}}W_i(\bar{k}_i+1,e_i+\sum_{j=1}^{l_i}\theta_{i,j}),W_i(\bar{k}_i+1,h_{v_i}) \Big\}\\
   &\le\tilde{\lambda}_i\max\Big\{\frac{\tilde{\lambda}_i^{l_i}}{\lambda_{W,i}^{l_i}}W_i(\bar{k}_i,e_i),\frac{\tilde{\lambda}_i^{l_i-1}}{\lambda_{W,i}^{l_i-1}}W_i(\bar{k}_i,e_i+\theta_{i,1}),\\
   &{\text{\quad\quad}}\dots,W_i(\bar{k}_i,e_i+\sum_{j=1}^{l_i}\theta_{i,j}),W_i(\bar{k}_i,e_i+\sum_{j=1}^{l_i}\theta_{i,j}) \Big\},\\
   &\le \tilde{\lambda}_i \tilde{W}_i(\bar{k}_i,l_i,\theta_i, e_i)
\end{split}
\end{equation}
where the first inequality is base on Items 2) and 3) in Assumption \ref{asm of protocol} and $h_{v_i}:=h_{v_i}(\bar{k}_i,e_i+\sum_{j=1}^{l}\theta_{i,j})$; while the second inequality is due to $\tilde{\lambda}_i>\lambda_i$ and $\lambda_{W,i}\ge1$. On the other hand, if the current instant is an updating instant, one has
\begin{equation}
\nonumber
  \begin{split}
  \tilde{W}_i&(\bar{k}_i,l_i-1,[\theta_{i,2}^{\rm{T}},\dots,\theta_{i,D_i+1}^{\rm{T}},0^{\rm{T}}]^{\rm{T}}, e_i+\theta_{i,1})\\
  &= \max\Big\{\frac{\tilde{\lambda}_i^{l_i-1}}{\lambda_{W,i}^{l_i-1}}W_i(\bar{k}_i,e_i+\theta_{i,1}),\dots,W_i(\bar{k}_i,e_i+\sum_{j=1}^{l_i}\theta_{i,j}) \Big\}\\
   &\le \tilde{W}_i(\bar{k}_i,l_i,\theta_i, e_i),
\end{split}
\end{equation}
where the equality uses the fact of $\theta_{i,j}=0$ with $j=l_i+1,\dots,D_i+1$. Thus, the inequalities in (\ref{jump of wi}) are proved.

Finally, consider the flow behavior characterized by (\ref{flow of wi}). For given $l_i$, suppose that $\tilde{W}_i(\bar{k}_i,l_i,\theta_i, e_i)=\frac{\tilde{\lambda}_i^{l_i-m}}{\lambda_{W,i}^{l_i-m}}W_i(\bar{k}_i,e_i+\sum_{j=1}^{m}\theta_{i,j})$ for some $m\in\{0,\dots, l_i\}$, that is, the value of $\frac{\tilde{\lambda}_i^{l_i-m}}{\lambda_{W,i}^{l_i-m}}W_i(\bar{k}_i,e_i+\sum_{j=1}^{m}\theta_{i,j})$ is maximum. Then, from Item 4) in Assumption \ref{asm of protocol}, Assumption \ref{asm of fvi}, and Lemma \ref{lem of clarke}, it follows that  for almost all $e_i\in\mathbb{R}^{n_{v,i}}$,
\begin{equation}
  \nonumber
  \begin{split}
    &\left\langle\frac{\partial \tilde{W}_i}{\partial e_i }, -f_{v_i}(q,w)\right\rangle\\
    =& \frac{\tilde{\lambda}_i^{l_i-m}}{\lambda_{W,i}^{l_i-m}} \left\langle\frac{\partial {W}_i(\bar{k}_i,e_i+\sum_{j=1}^{m}\theta_{i,j})}{\partial (e_i+\sum_{j=1}^{m}\theta_{i,j})}, -f_{v_i}(q,w)\right\rangle \\
    \le& \frac{\tilde{\lambda}_i^{l_i-m}}{\lambda_{W,i}^{l_i-m}}\left(H_i(x,e,w)+M_{e,i}\left\|e_i\right\|\right) \\
    \le& H_i(x,e,w)+\frac{\tilde{\lambda}_i^{l_i-m}M_{e,i}}{\lambda_{W,i}^{l_i-m}\underline{\beta}_{W_i}}W_i(\bar{k}_i, e_i),
  \end{split}
\end{equation}
where the last inequality is based on $\frac{\tilde{\lambda}_i}{\lambda_{W,i}}<1$ and Item 1) in Assumption \ref{asm of protocol}. From the definition in (\ref{def of tilde wi}), it follows
\[\frac{\tilde{\lambda}_i^{l_i}}{\lambda_{W,i}^{l_i}}W_i(\bar{k}_i,e_i)\le\frac{\tilde{\lambda}_i^{l_i-m}}{\lambda_{W,i}^{l_i-m}}W_i(\bar{k}_i,e_i+\sum_{j=1}^{m}\theta_{i,j})=\tilde{W}_i, \]
which implies
\begin{equation}
  \nonumber
  \begin{split}
    \left\langle\frac{\partial \tilde{W}_i}{\partial e_i}, -f_{v_i}(q,w)\right\rangle\le &H_i(x,e,w)+\frac{\lambda_{W,i}^{m}M_{e,i}}{\tilde{\lambda}_i^{m}\underline{\beta}_{W_i}}\tilde{W}_i,\\
    \le & H_i(x,e,w)+\frac{\lambda_{W,i}^{l_i}M_{e,i}}{\tilde{\lambda}_i^{l_i}\underline{\beta}_{W_i}}\tilde{W}_i,
  \end{split}
\end{equation}
which completes the proof.
\end{IEEEproof}

\begin{IEEEproof}[Proof of Proposition \ref{pps of V}]
  From Proposition \ref{pps of we} and Items 1) and 3) in Assumption \ref{asm of V orignal}, one can directly prove (\ref{positive of V}) and (\ref{flow of deltai}) since $\tilde{V}=V(x)$ and $\tilde{\delta_i}=\delta_i$. For the flow behavior of $\tilde{V}=V(x)$, Item 2) in Assumption \ref{asm of V orignal} implies
  \begin{equation}
  \nonumber
  \begin{split}
     \left\langle \nabla \tilde{V}(x), f(q,w)\right\rangle\le& -\alpha_{\tilde{V}}(\left\|x\right\|)+\sum_{i=1}^N\Big((\gamma_{i}^2-\bar{\varepsilon}_i)W_i^2\\
     &-\tilde{\delta}_i(v_i)-(1+\epsilon_{l_i,i})H_{l_i,i}^2(x,e,w)\\
     &-\hat{\delta}_i(\hat{v}_i)-(1+\epsilon_{l_i,i})\tilde{J}_i(x,e,w)\Big)\\
     &+\alpha_w(\left\|w\right\|),
  \end{split}
\end{equation}
where we utilize the equalities in (\ref{notations in pps V}). For each given $l_i$, from the definition of $\tilde{W}_i(\bar{k}_i,l_i,\theta_i, e_i)$ in (\ref{def of tilde wi}), it follows
\[  W_i(\bar{k}_i,e_i)\le \frac{\lambda_{W,i}^{l_i}\tilde{W}_i(\bar{k}_i,l_i,\theta_i, e_i)}{\tilde{\lambda}_i^{l_i}}.\]
Thus, we have
  \begin{equation}
  \nonumber
  \begin{split}
     \left\langle \nabla \tilde{V}(x), f(q,w)\right\rangle\le& -\alpha_{\tilde{V}}(\left\|x\right\|)+\sum_{i=1}^N\Big(\frac{\lambda_{W,i}^{2l_i}(\gamma_{i}^2-\bar{\varepsilon}_i)}{\tilde{\lambda}^{2l_i}}\tilde{W}_i^2\\
     &-\tilde{\delta}_i(v_i)-(1+\epsilon_{l_i,i})H_{l_i,i}^2(x,e,w)\\
     &-\hat{\delta}_i(\hat{v}_i)-(1+\epsilon_{l_i,i})\tilde{J}_i(x,e,w)\Big)\\
     &+\alpha_w(\left\|w\right\|),
  \end{split}
\end{equation}
which completes the proof.
\end{IEEEproof}

\begin{IEEEproof}[Proof of Proposition \ref{pps of upper bound of l}]
For each communication channel $\mathcal{C}_i,i\in\bar{N}$, due to Assumption \ref{asm of large delay}, at every sampling instant $s_j^i$, the transmission at or before $s_{j-D_i-1}^i$  must have been updated. In the small-delay case, we have $l_i(s_j^i)=\hat{l}_i(s_j^i)=D_i=0$ for all $j\in\mathbb{Z}_{\ge0}$, since all the updates have to be executed before the next sampling instant. Otherwise, to study the upper bound of $l_i(t)$ at $t=s_{j+1}^i$, we only need to count the number of transmissions within $\{s_{j-D_i+1}^i,\dots,s_j^i\}$ whose cardinality is $D_i$. Therefore, the proof is completed from the definition of $\hat{l}_i(s_j^{j+})$.
\end{IEEEproof}

\end{document}